%% file: main.tex
\documentclass[12pt]{article}  
\usepackage[a4paper, total={16cm, 23cm}]{geometry}

\input{pkgssetup.tex}

\title{Data-driven optimal approximation on Hardy spaces in simply connected domains}
		\author{Alessandro Borghi\thanks{borghi@tu-berlin.de, Institut für Mathematik, Technische Universität Berlin, Straße des 17.
Juni 136, Berlin, 10623, Germany}  and Tobias Breiten\thanks{tobias.breiten@tu-berlin.de, Institut für Mathematik, Technische Universität Berlin, Straße des 17.
Juni 136, Berlin, 10623, Germany}} 
		\date{}

\begin{document}

\maketitle

\section*{Abstract}
		We consider optimal interpolation of functions analytic in simply connected domains in the complex plane.
		By choosing a specific structure for the approximant, we show that the resulting first order optimality conditions can be interpreted as optimal $\mathcal{H}_2$ interpolation conditions for discrete-time dynamical systems. 
        Connections to the implicit Euler method, the midpoint method, and backward differentiation methods are also established.
		A data-driven algorithm is developed to compute a (locally) optimal approximant.  
        Our method is tested on three numerical experiments.
\newpage


\section{Introduction}
In this article we consider the optimization problem 
\begin{equation}\label{eq:generaloptprob}
    \min_{G_r\in\mathcal{X}}\|G-G_r\|_{\mathcal{X}},
\end{equation}
for a given a function $G$ and some Hardy space $\mathcal{X}$ over a fixed, arbitrary simply connected domain in $\mathbb C$.
Our main interest for studying \eqref{eq:generaloptprob} is its relevance in the context of optimal model order reduction. In particular, we discuss the construction of approximants $G_r$ of the form 
\begin{equation}\label{eq:weirdapprox}
    \begin{aligned}
        G_r(\cdot) &= \left(\bfc_r^*\left(\varphi^{-1}(\cdot)\bfI_r - \bfA_r\right)^{-1}\bfb_r + d_r\right)[\varphi'(\varphi^{-1}(\cdot))]^{-1/2},  \\
        \text{ or }G_r(\cdot) &= \bfc_r^*\left(\varphi^{-1}(\cdot)\bfI_r - \bfA_r\right)^{-1}\bfb_r + d_r,     
    \end{aligned}
\end{equation}    
where $\bfA_r\in\C^{r\times r}$, $d_r\in\C$, $\bfb_r,\bfc_r\in\C^{r}$, and $\varphi$ is a specific conformal map defined in more detail later on.
For our framework, we draw upon and generalize existing literature on interpolatory $\hardy$-optimal model order reduction, see, e.g., \cite{AnBeGu20,Ant05}.
Well known and studied examples are linear time invariant (LTI) systems of the form 
\begin{equation}\label{eq:sscont}
    \begin{cases}
        x(t) &= \bfA x(t) + \bfb u(t), \quad x(0) = 0,\\
        y(t) &= \bfc^* x(t) + du(t), 
    \end{cases}
\end{equation}
where $\bfA\in\C^{n\times n}$, $d\in\C$, and $\bfb,\bfc\in\C^{n}$. 
The system in \eqref{eq:sscont} is referred to as the full order model (FOM). 
Here, $n$ is large, making it computationally challenging to obtain the dynamics of the system. 
The goal is to compute a reduced order model (ROM) of the form
\begin{equation*}
    \begin{cases}
        \widehat{x}(t) &= \bfA_r \widehat{x}(t) + \bfb_r u(t), \quad \widehat{x}(0) = 0,\\
        \widehat{y}(t) &= \bfc_r^* \widehat{x}(t) + d_ru(t), 
    \end{cases}
\end{equation*}
with $r\ll n$, such that $y(t)\approx\widehat{y}(t)$ for a wide range of inputs $u$.
To this end, $\hardy$-optimal model order reduction solves \eqref{eq:generaloptprob} for the Hardy space $\mathcal{X}=\hardy(\C_+)$ of functions analytic in the right half plane.
This translates into minimizing the error norm
\[
\|G-G_r\|_{\hardy(\C_+)}^2 = \frac{1}{2\pi}\int_{-\infty}^{\infty} |G(\iunit\omega) - G_r(\iunit\omega)|^2\mathrm{d}\omega,
\]
where $G$ and $G_r$ are the transfer functions of the FOM and ROM, i.e.,
\[
G(\cdot) = \bfc^*(\cdot\bfI - \bfA)^{-1}\bfb + d, \textnormal{ and } G_r(\cdot) = \bfc_r^*(\cdot\bfI_r - \bfA_r)^{-1}\bfb_r + d_r. 
\]
Of particular importance is the $\mathcal{L}_\infty$ error bound that the framework provides (see, e.g., \cite[Section 2.1.1]{AnBeGu20}) 
\begin{equation}\label{eq:Linftyboundcont}
\|y-\widehat{y}\|_{\mathcal{L}_\infty} \leq\|G-G_r\|_{\hardy(\C_+)}\|u\|_{\mathcal{L}_2},
\end{equation}
and which motivates computing a ROM minimizing the $\hardy(\C_+)$ error norm. 
First order optimality conditions for the $\hardy$-optimal model reduction of LTI systems were derived both in continuous \cite{MeLu67,Wil70} and discrete-time \cite{BunKuVoWi10}.
This resulted in  numerical algorithms for the computation of ROM minimizing the $\hardy$ error norm.
In particular, the iterative rational Krylov algorithm (IRKA) proposed in \cite{GuAnBe08} and MIRIAm from \cite{BunKuVoWi10} were designed for continuous and discrete-time systems, respectively.
While IRKA is a projection based algorithm, a data-driven version, transfer function IRKA or TF-IRKA, based on Loewner matrices was also proposed in \cite{BeGu12}.
Several extensions of IRKA for second order systems \cite{Wya12}, bilinear systems \cite{BeBr12}, linear systems with quadratic output \cite{ReDuGoGu24, ReGoDuGu25}, port-Hamiltonian systems \cite{GuPoBeSc12}, and differential algebraic equations \cite{GuStWy13} were developed over the years.
A unifying framework for $\hardy$-optimal model reduction of structured time invariant systems has lately been proposed in \cite{MlBeGu25}, earlier related results are found in \cite{MlGu23,MlGu23b}. 
Additionally, in \cite{BorBre23} an optimal model reduction framework was introduced for LTI systems with poles in domains in the complex plane that are not necessarily the unit disk or the left half plane. 

In this work, we take inspiration from \cite{BorBre23} and build upon the $\hardy$-optimal model reduction framework to solve problems as in \eqref{eq:generaloptprob}. 
In particular  we consider  $\mathcal{X}$ as a Hardy space of functions analytic in simply connected domains.
More precisely, we rely on \cite{Dur70} which generalizes the classical Hardy space $\hardy$ in two different ways referred to as $\etwod$ and $\durenhardy$.
Our main contributions are the following:
\begin{enumerate}
    \item We show that by defining $G_r$ as in \eqref{eq:weirdapprox}, we can recover explicit optimality conditions in $\etwod$ and $\durenhardy$. 
    \item We connect our framework to the computation of reduced models of discretized LTI systems and, similar to \eqref{eq:Linftyboundcont}, prove that there exists an $\ell_\infty$ error bound based on the $\durenhardy$ error norm between the respective transfer functions. 
    \item Based on the resulting optimality conditions, we develop a TF-IRKA based algorithm for the computation of (locally) optimal approximants in $\etwod$ and $\durenhardy$.
\end{enumerate}
The rest of the paper is organized as follows. 
In \cref{sec:optMORh2} we briefly recapitulate the $\hardy$-optimal model order reduction framework for discrete-time systems along with some extensions. 
In \cref{sec:durenspace} we introduce the concept of Hardy spaces of functions analytic in general domains. In \cref{sec:3} we develop an optimal $\etwod$ and $\durenhardy$ model reduction framework and derive first order necessary optimality conditions. 
Additionally, we design a TF-IRKA based algorithm and provide a connection with model order reduction of discretized LTI systems. 
We conclude with \cref{sec:numerics} where we test our algorithm on three numerical experiments. 

\subsection{Notation}
Throughout the paper we denote the complex plane by $\C$. 
If $\A$ is a simply connected domain, $\partial\A$ is its boundary, $\bar{\A}=\{\A\cup\partial\A\}$ is its closure, $\A^{\mathsf{c}}$ is its complement, and $\bar{\A}^\mathsf{c}=\{\C\backslash\bar{\A}\}$ is its exterior. 
In addition, the symbols $\C_-$, $\C_+$, and $\unitdisk$, stand for the left-half complex plane, the right-half complex plane, and the unit disk. 
$\R$ and $\iunit\R$ indicate the real and imaginary numbers, respectively. 
For a single-variable complex valued differentiable bijective function $f$ we indicate its complex derivative as $f'$, and its inverse as $f^{-1}$.
The composition of two functions $f$ and $g$ is indicated by $f(g(\cdot))$ or $f\circ g$.
The symbol $\overline{(\cdot)}$ indicates the complex conjugate of a scalar and the absolute value of a complex number $z$ is denoted by $|z|=\sqrt{z\overline{z}}$. 
For the transpose and conjugate transpose of a vector or matrix we use the symbols $(\cdot)^\top$ and $(\cdot)^{*}$, respectively. 
When declaring a Hardy space of functions we explicitly write the domain where its elements are analytic.
For example $\hardy(\A)$ is the Hardy space of functions analytic in $\A$.


\section{Preliminaries}
In this section we recall the concept of interpolatory model order reduction, specifically the problem of optimal model reduction for linear discrete-time invariant systems and its resulting optimal interpolation conditions.
We then discuss the latest state-of-the-art methods and discoveries on optimal model reduction in Hardy spaces.

\subsection{Optimal model order reduction on Hardy spaces}\label{sec:optMORh2}
Optimal model order reduction deals with the computation of solutions to problems of the form \eqref{eq:generaloptprob}. Let us for example consider a discrete-time LTI dynamical system of the form 
\begin{equation}\label{eq:LTIdynsys}
    \begin{cases}
        x_{k+1} &= \bfA x_{k} + \bfb u_{k}, \quad x_0 = 0,\\
        y_{k} &= \bfc^* x_k + du_k, 
    \end{cases} 
\end{equation}
where $\bfA\in\C^{n\times n}$, $d\in\C$, and $\bfb,\bfc\in\C^{n}$. In addition, $x_k\in\C^n$, $u_{k}\in\C$, and $y_k\in\C$ are known as the states, inputs, and outputs of the system, respectively, at step $k$. Its frequency domain transfer function
\begin{equation}\label{eq:LTItf}
    G(\cdot) = \bfc^* (\cdot \bfI - \bfA)^{-1}\bfb + d,
\end{equation}
resulting from applying the $\mathcal{Z}$-transform, is a pivotal tool in the reduction of \eqref{eq:LTIdynsys}.
If the system \eqref{eq:LTIdynsys} is asymptotically stable all eigenvalues of $\bfA$ and, consequently, all poles of $G$ lie inside the unit disk $\unitdisk$.
Since $G$ is a rational function, it is therefore analytic in $\bar{\unitdisk}^\mathsf{c}$. In this case, a canonical choice for $\mathcal{X}$ is the Hardy space
\begin{equation}\label{eq:H2Dc}
\hardy(\bar{\unitdisk}^{\mathsf{c}}):=\left\{G\colon\bar{\unitdisk}^{\mathsf{c}}\rightarrow\C \textnormal{ analytic }\Big|\sup_{r>1}\int_{0}^{2\pi}|G(re^{\iunit\theta})|^2\mathrm{d}\theta<\infty\right\},
\end{equation}
see, e.g., \cite[Chapter 5.1.3]{Ant05}.
Given $G\in\hardy(\bar{\unitdisk}^\mathsf{c})$, consider a reduced order model
\[
G_r(\cdot) = \bfc_r^*(\cdot\bfI_r -\bfA_r)^{-1}\bfb_r +d_r \in \hardy(\bar{\unitdisk}^\mathsf{c}),
\]
where $\bfA_r\in\C^{r\times r}$, $d_r\in\C$, and $\bfb_r,\bfc_r\in\C^{r}$, with $r\ll n$ and poles in $\unitdisk$. The optimization problem \eqref{eq:generaloptprob}  now reads 
\begin{equation}\label{eq:optprobH2}
\min_{G_r\in\hardy(\bar{\unitdisk}^\mathsf{c})}\|G-G_r\|_{\hardy(\bar{\unitdisk}^\mathsf{c})}.
\end{equation}
In \cite{GuAnBe08,BunKuVoWi10} it was proved that if $G_r$ is a local minimizer of \eqref{eq:optprobH2} with $d=d_r$ then the following Hermite interpolation conditions hold
\begin{equation}\label{eq:intoptH2}
    G(1/\overline{\lambda_j}) = G_r(1/\overline{\lambda_j}) \textnormal{ and } G'(1/\overline{\lambda_j}) = G_r'(1/\overline{\lambda_j}) \quad \text{for }j=1,\dots,r,
\end{equation}
where $\lambda_j\in\unitdisk$, $j=1,\dots,r$, are the poles of $G_r$.

Similarly, for asymptotically stable linear continuous-time systems as in \eqref{eq:sscont}, one typically considers the Hardy space of functions analytic in the right-half plane
\begin{equation}\label{eq:H2Cplus}
\hardy(\C_+):=\left\{G\colon\C_+\rightarrow\C \textnormal{ analytic }\Big|\sup_{x>0}\int_{-\infty}^{\infty}|G(x+\iunit y)|^2\mathrm{d}y<\infty\right\}.
\end{equation}
For the corresponding optimization problem \eqref{eq:generaloptprob}, we have the following optimality conditions (see, e.g.,  \cite{MeLu67,GuAnBe08,GuStWy13})
\begin{equation}\label{eq:intoptH2Cplus}
    G(-\overline{\lambda_j}) = G_r(-\overline{\lambda_j}), \textnormal{ and } G'(-\overline{\lambda_j}) = G_r'(-\overline{\lambda_j}) \quad \text{for }j=1,\dots,r,
\end{equation}
with $\lambda_j\in\C_-$ and $d=d_r$. 
Recently in \cite{MlBeGu25} interpolatory $\hardy(\C_+)$ optimality conditions were obtained for a wider range of dynamical systems.
More in detail, the authors consider diagonally structured time invariant systems which include time-delay and second order systems. 
While \cite{MlBeGu25} specifically focused on dynamical systems in the $\hardy(\C_+)$ space, more general necessary optimality conditions were found in previous works \cite{MlGu23,MlGu23b} for reduced order models that minimize an $\mathcal{L}_2$ error norm.

For cases where the transfer function does neither belong to $\hardy(\bar{\unitdisk}^\mathsf{c})$ nor $\hardy(\C_+)$ the conditions \eqref{eq:intoptH2} and \eqref{eq:intoptH2Cplus} are not suitable anymore.
As a remedy, \cite{Kub08} introduced the $h_{2,\alpha}$ and $\mathcal{H}_{2,\alpha}$ frameworks which can handle LTI systems with poles inside a disk of radius $\alpha$ and in an $\alpha$-shifted left-half plane, respectively. In particular, optimality conditions similar to \eqref{eq:intoptH2} and \eqref{eq:intoptH2Cplus} can be derived for a wider range of transfer functions, including unstable systems. 

An optimal model reduction framework for functions analytic in general domains was also proposed in \cite{BorBre23}.
Here, the concept of Hardy spaces in general domains was used to design an IRKA based algorithm. For an appropriate generalization of the conditions in \eqref{eq:intoptH2} and \eqref{eq:intoptH2Cplus}, conformal maps and their characterization, which we briefly recall below, were used.
\begin{theorem}(\cite{Weg12}, Theorem 6.1.2)\label{th:confmap}
    Suppose $\mathbb{X},\mathbb{Y}\subset\mathbb{C}$ are open sets and let $\varphi\colon\mathbb{X}\rightarrow\mathbb{Y}$ be Fr\'echet differentiable as a function of two real variables. The mapping $\varphi$ is conformal in $\mathbb{X}$ if and only if it is analytic in $\mathbb{X}$ and $\varphi'(z_0)\neq 0$ for every $z_0\in\mathbb{X}$.
\end{theorem}
Under suitable assumptions on $\varphi$, the transfer function $G$ with poles in $\mathbb{P}\subset\C$ can be mapped to a function $\calA_G(\cdot):=G(\varphi(\cdot))[\varphi(\cdot)]^{1/2}$ analytic in $\C_+$.
The considered space of functions is
\[
\hardy(\bar{\mathbb{P}}^{\mathsf{c}}):=\left\{G\colon\bar{\mathbb{P}}^{\mathsf{c}}\rightarrow\C \text{ analytic }\Big| \|G\|_{\hardy(\bar{\mathbb{P}}^{\mathsf{c}})}<\infty\right\},
\]
with 
\[
    \|G\|_{\hardy(\bar{\mathbb{P}}^{\mathsf{c}})} = \left(\frac{1}{2\pi}\int_{-\infty}^{\infty} |\calA_G(\iunit y)|^2\mathrm{d}y\right)^{1/2}.
\]
For $G_r\in\hardy(\bar{\mathbb{P}}^{\mathsf{c}})$ being a minimizer of $\|G-G_r\|_{\hardy(\bar{\mathbb{P}}^{\mathsf{c}})}$, the following interpolation conditions hold (see \cite[Corollary 3]{BorBre23})
\begin{equation}\label{eq:intoptH2general}
    G\left(\varphi(-\overline{\varphi^{-1}(\lambda_j)})\right) = G_r\left(\varphi(-\overline{\varphi^{-1}(\lambda_j)})\right), \textnormal{ and } G'\left(\varphi(-\overline{\varphi^{-1}(\lambda_j)})\right) = G_r'\left(\varphi(-\overline{\varphi^{-1}(\lambda_j)})\right)
\end{equation}
for $j=1,\dots,r$.
Even though \eqref{eq:intoptH2general} provides a basis for generalizations of IRKA, the derivations in \cite{BorBre23} were carried out under a rather restrictive set of assumptions on $\varphi$. For this reason, in what follows we discuss an alternative optimization framework which relies on a particular structure of the sought approximations $G_r$. As a result, we obtain optimality conditions similar to \eqref{eq:intoptH2}.

\section{The $\etwod$ and $\durenhardy$ spaces} \label{sec:durenspace}
In this section, we collect some well-known results about Hardy spaces in general domains for which we closely follow the exposition in \cite{Dur70}.
\begin{definition}[$\hardy$ space]\label{def:H2}
    Let $G,F$ be analytic functions in the unit disk $\unitdisk$. Then $G$ is of class $\hardy$ if 
    \begin{equation*}
        \sup_{r\nearrow 1} \frac{1}{2\pi}\int_0^{2\pi}|G(re^{\iunit\theta})|^2\mathrm{d}\theta<\infty.
    \end{equation*}
    We define a space as 
    \begin{equation*}
        \hardy:=\left\{G\colon \unitdisk\rightarrow \C \textnormal{ analytic}\Big| \sup_{r\nearrow 1} \frac{1}{2\pi}\int_0^{2\pi}|G(re^{\iunit\theta})|^2\mathrm{d}\theta<\infty \right\}.
    \end{equation*}
    For $G,F\in\hardy$ we have the $\hardy$ inner product
    \[
    \langle G, F\rangle_{\hardy} = \frac{1}{2\pi}\int_0^{2\pi}G(e^{\iunit\theta})\overline{F(e^{\iunit\theta})}\mathrm{d}\theta,
    \]
    and induced norm 
    \begin{equation*}
        \|G\|_{\hardy}= \sqrt{\langle G,G\rangle_{\hardy}} = \left(\frac{1}{2\pi}\int_0^{2\pi}|G(e^{\iunit\theta})|^2\mathrm{d}\theta\right)^{1/2}.
    \end{equation*}
\end{definition}
With regard to \cref{sec:optMORh2} and the Hardy space $\hardy(\bar{\unitdisk}^\mathsf{c})$, let us emphasize that a transfer function $G\in\hardy(\bar{\unitdisk}^\mathsf{c})$ in \eqref{eq:LTItf} can be related to a function in $\hardy$ via the transformation $z\rightarrow z^{-1}$ which guarantees $G((\cdot)^{-1})$ to be analytic in the unit disk (see \cite[Section 3]{ChuChe97}).

Based on the characterization in \cite[Section 10.1]{Dur70}, below we provide definitions to two possible extensions of Hardy spaces for functions analytic in general domains $\A$.
\begin{definition}[$\etwod$ space]\label{def:hardyspacegeneral1}
    Let $G\colon \A\rightarrow \C$ and $F\colon \A\rightarrow \C$ be analytic and $\varphi\colon\unitdisk\rightarrow \A$ be conformal with $\A$ a simply connected domain. 
    Let 
    \begin{equation}\label{eq:AG}
    \mathcal{A}_G(\cdot):=G(\varphi(\cdot)) [\varphi'(\cdot)]^{1/2},
    \end{equation}
    we then define the Hardy space 
    \begin{align*}
        \etwod := \left\{G\colon\A\rightarrow \C \textnormal{ analytic } \big| \sup_{r\nearrow 1}\frac{1}{2\pi}\int_{0}^{2\pi}|\mathcal{A}_G(re^{\iunit\theta})|^2\mathrm{d}\theta<\infty\right\},   
    \end{align*} 
    with inner product $\langle G,F\rangle_{\etwod} = \langle \mathcal{A}_G, \mathcal{A}_F \rangle_{\hardy}$,
    and induced norm $\|G\|_{\etwod} = \|\mathcal{A}_G\|_{\hardy}$.
\end{definition}
\begin{definition}[$\durenhardy$ space]\label{def:hardyspacegeneral2}
    Let $G\colon \A\rightarrow \C$ and $F\colon \A\rightarrow \C$ be analytic and $\varphi\colon\unitdisk\rightarrow \A$ be conformal with $\A$ a simply connected domain. 
    We then define the Hardy space
    \begin{align*}
        \durenhardy := \left\{G\colon\A\rightarrow \C \textnormal{ analytic } \big| \sup_{r\nearrow 1}\frac{1}{2\pi}\int_{0}^{2\pi}|G\circ\varphi(re^{\iunit\theta})|^2\mathrm{d}\theta<\infty\right\},
    \end{align*}
    with inner product $
    \langle G,F\rangle_{\durenhardy} = \langle G\circ\varphi, F\circ\varphi\rangle_{\hardy},  
    $
    and induced norm $\|G\|_{\durenhardy} = \|G\circ\varphi\|_{\hardy}$.
\end{definition}

For more details on Hardy spaces in general domains we refer the reader to \cite[Chapter 10]{Dur70}.
Throughout this manuscript we make the following assumption.
\begin{assumption}\label{assumption:conformalmap}
    Let $\varphi\colon\unitdisk\rightarrow \A$ be a bijective conformal map.
\end{assumption}

Let us take a closer look at the evaluation of the $\etwod$ and $\durenhardy$ inner products when one of the involved functions has a particular structure.
For the $\etwod$ space we consider 
\begin{equation}\label{eq:rho}
    \rho(\cdot)= \frac{1}{(\varphi^{-1}(\cdot) - \lambda)[\varphi'(\varphi^{-1}(\cdot))]^{1/2}} \quad \text{with }\lambda\in\C \text{ and } |\lambda|>1   . 
\end{equation}
It is worth mentioning that for $|\lambda|>1$ we have that $\rho\in \etwod$ as
\[
\mathcal{A}_\rho=\rho(\varphi(\cdot))[\varphi'(\cdot)]^{1/2} = \frac{1}{\cdot - \lambda}\in \hardy, 
\]
where we used \eqref{eq:AG}.
For $G\in \durenhardy$, instead, we consider the function
\begin{equation}\label{eq:rho2}
\rho(\cdot) = \frac{1}{\varphi^{-1}(\cdot) - \lambda}, \quad \text{with } \lambda\in\C \text{ and }|\lambda|>1.
\end{equation}
We then have the following result.
\begin{lemma}\label{lemma:E2DH2Dinnerprodeval}
    Let $G\in\etwod$ and $\rho\in\etwod$ as in \eqref{eq:rho}. Then 
    \begin{equation}\label{eq:evalinnerproductE2D}
        \left\langle G,\rho\right\rangle_{\etwod} = -\frac{1}{\overline{\lambda}}\calA_G\left(\frac{1}{\overline{\lambda}}\right),
    \end{equation}
    with $\calA_G$ as in \eqref{eq:AG}.
    For $G\in\durenhardy$ and $\rho\in\durenhardy$ as in \eqref{eq:rho2}, instead, we have
    \begin{equation}\label{eq:evalH2Dinprod}
			\left\langle G,\rho\right\rangle_{\durenhardy} = -\frac{1}{\overline{\lambda}}G\circ\varphi\left(\frac{1}{\overline{\lambda}}\right).
	\end{equation}
\end{lemma}
\begin{proof}
By transforming the $\etwod$ inner product into an $\hardy$ inner product (see \cref{def:H2}) and applying a change of coordinates $z=e^{\iunit\theta}$ followed by the residue theorem we get
\begin{equation*}
    \begin{aligned}
        \left\langle G,\rho\right\rangle_{\etwod} &= \left\langle \calA_G, \calA_\rho \right\rangle_{\hardy} = \left\langle \calA_G, \frac{1}{\cdot - \lambda} \right\rangle_{\hardy} = -\frac{1}{\overline{\lambda}}\calA_G\left(\frac{1}{\overline{\lambda}}\right),
    \end{aligned}
\end{equation*} 
proving \eqref{eq:evalinnerproductE2D}.
A similar result can be also found in \cite[Section 6.5.1]{Nik19}.
The proof for \eqref{eq:evalH2Dinprod} follows the same steps.
\end{proof}

We will show in the next section that for both spaces it is possible to design an optimal approximation framework.
Here, by choosing weighted sums of $\rho$ as approximants we can retrieve first order optimality conditions similar to \eqref{eq:intoptH2}.

\section{The $\etwod$ and $\durenhardy$ approximation frameworks} \label{sec:3}
In this section we design an optimal approximation framework for functions in $\etwod$ and $\durenhardy$. 
The two frameworks are divided due to the different assumptions on the approximant, however, the two spaces coincide if $|\varphi'(\cdot)|$ is bounded from above and below (see \cite[Theorem 10.2]{Dur70}).
While similar work has been done for optimal approximation on Hardy spaces in general domains \cite{BorBre23,Kub08}, here we do not assume the approximant to have a rational structure. 
The rational structure is visible once the approximant is recasted into the $\hardy$ space through the conformal map $\varphi$.

\subsection{Optimal interpolation conditons in the $\etwod$ and $\durenhardy$ frameworks}\label{sec:approx}

We start with the $\etwod$ optimal framework as $\durenhardy$ can be easily derived from the former.
Motivated by the results from  \cref{sec:durenspace}, consider an approximant $G_r$ of the form
\begin{equation}\label{eq:Gr}
    G_r(\cdot) = \left(\sum_{j=1}^r \frac{\phi_j}{(\varphi^{-1}(\cdot) - \lambda_j)} +d_r\right)\frac{1}{[\varphi'(\varphi^{-1}(\cdot))]^{1/2}},
\end{equation}
with $\lambda_j,d_r,\phi_j\in\C$, and  $|\lambda_j|>1,\forall j=1,\dots ,r$. Note that $G_r\in \etwod$.
Here $\calA_{G_r}$ relates to the transfer function of a linear time invariant system  since we have
\[
G_r(\cdot) = \left(\bfc_r^*(\varphi^{-1}(\cdot)\bfI_r-\mathbf{\Lambda}_r)^{-1}\bfb_r+d_r\right)[\varphi'(\varphi^{-1}(\cdot))]^{-1/2},
\]
with $\bfc_r,\bfb_r\in\C^r$ such that $(\bfc_r^*)_i(\bfb_r)_i=\phi_i$, and $\mathbf{\Lambda}_r=\texttt{diag}(\lambda_1,\dots,\lambda_r)$. Hence, we have that
\begin{equation}\label{eq:AGr}
    \calA_{G_r}(\cdot) =G_r(\varphi(\cdot))[\varphi'(\cdot)]^{1/2} = \bfc_r^*(\cdot\bfI_r-\mathbf{\Lambda}_r)^{-1}\bfb_r+d_r \in \hardy.
\end{equation}

The main objective of the optimal $\etwod$ approximation framework is to compute $G_r\in\etwod$  that minimizes the cost function $\|G-G_r\|_{\etwod}$.
We rewrite the optimization problem with respect to the coefficients of $G_r$ as follows
\begin{equation}\label{eq:optE2D}
    \min_{\substack{\lambda_j,d_r,\phi_j\in\C \\ |\lambda_j|>1 \\ j=1,\dots,r}}\|G-G_r\|_{\etwod}.
\end{equation}
In the next theorem we derive first order optimality conditions for \eqref{eq:optE2D}.
\begin{theorem}\label{th:E2Dintcond}
    Let $G_r\in \etwod$ in \eqref{eq:Gr} be a local minimizer of \eqref{eq:optE2D}. We then have the following interpolation conditions
    \begin{equation}\label{eq:optintE2D}
        \begin{aligned}
            G\circ\varphi\left(0\right) &= G_r\circ\varphi\left(0\right), \\ G\circ\varphi\left(\frac{1}{\overline{\lambda_p}}\right) &= G_r\circ\varphi\left(\frac{1}{\overline{\lambda_p}}\right) \\ \text{ and } G'\circ\varphi\left(\frac{1}{\overline{\lambda_p}}\right) &= G_r'\circ\varphi\left(\frac{1}{\overline{\lambda_p}}\right), 
        \end{aligned}
    \end{equation}
    for $p=1,\dots,r$.
\end{theorem}
\begin{proof}
    The proof uses standard arguments and follows along the lines of similar proofs in, e.g., \cite{AnBeGuAn13,AnBeGu20,BorBre23,GuAnBe08}. For a self-contained presentation, we provide the proof in \cref{proof:e2dintcond}.
\end{proof}
It is possible to see the resemblance between \eqref{eq:optintE2D} and the $\hardy(\bar{\unitdisk}^{\mathsf{c}})$ optimality conditions in \eqref{eq:intoptH2}.
This is due to the close connection between the $\etwod$ and the $\hardy$ spaces.
While the $\hardy$ and $\hardy(\bar{\unitdisk}^{\mathsf{c}})$ spaces are different, the resulting interpolation conditions are the same, i.e., a reflection of the reduced model poles with respect to the unit circle. 

\begin{remark} \label{remark:differenceborbre}
    Let us comment on the differences between the  $\etwod$ optimal framework and \cite{BorBre23}.
    Firstly, while in \cite{BorBre23} the approximant is defined as the transfer function of an LTI system, $G_r$ in \eqref{eq:Gr} is generally not a rational function.
    Additionally, the framework in \cite{BorBre23} limits the choice of the domain $\A$ due to the restrictive assumptions on the conformal maps.
    In this manuscript we only assume $\varphi$ to be conformal and bijective, see \cref{assumption:conformalmap}.
    Due to $\A$ being a simply connected domain, the Riemann mapping theorem \cite[Theorem 2.15]{Kyt19} asserts the existence of a bijective conformal map from $\A$ to $\unitdisk$.
    This provides more flexibility in the choice of $\A$ and $\varphi$.
\end{remark}

For the derivation of the optimality conditions for the $\durenhardy$ framework, we consider
\begin{equation}\label{eq:Gr2}
    G_r(\cdot) = \sum_{j=1}^r \frac{\phi_j}{\varphi^{-1}(\cdot) - \lambda_j} + d_r
\end{equation}
with $d_r,\phi_j,\lambda_j\in\C$, and $|\lambda_j|>1$, $\forall j=1,\dots,r$. 
Similar to the $\etwod$ framework, we want to compute $G_r\in\durenhardy$ such that it solves 
\begin{equation}\label{eq:H2Dopt}
	\min_{\substack{\lambda_j,d_r,\phi_j\in\C \\ |\lambda_j|>1 \\ j = 1,\dots,r}}\|G-G_r\|_{\durenhardy}.
\end{equation}
Similar to \cref{th:E2Dintcond}, for a minimizer $G_r$ of \eqref{eq:H2Dopt} we have the following interpolation conditions.
\begin{theorem}\label{theorem:H2Dintcond}
    Let $G_r\in \durenhardy$ in \eqref{eq:Gr2} be a local minimizer of \eqref{eq:H2Dopt}. We then have the following interpolation conditions
    \begin{equation}\label{eq:optintH2D}
        \begin{aligned}
            G\circ\varphi\left(0\right) &= G_r\circ\varphi\left(0\right) \\ G\circ\varphi\left(\frac{1}{\overline{\lambda_p}}\right) &= G_r\circ\varphi\left(\frac{1}{\overline{\lambda_p}}\right), \\ \text{ and } G'\circ\varphi\left(\frac{1}{\overline{\lambda_p}}\right) &= G_r'\circ\varphi\left(\frac{1}{\overline{\lambda_p}}\right), 
        \end{aligned}
    \end{equation}
    for $p=1,\dots,r$.
\end{theorem}
\begin{proof}
	The proof follows the same steps used for \cref{th:E2Dintcond}.
    The main difference is in the choice of the approximant, which in this case is \eqref{eq:Gr2}.
\end{proof}

It is interesting to see that for specific choices of $\varphi^{-1}$ the approximant $G_r$ can be viewed as the transfer function of a diagonally structured time invariant system 
\[
G_r(z) = \mathcal{C}_r(z)\mathcal{K}_r(z)^{-1}\mathcal{B}_r(z) = \sum_{j=1}^{r}\frac{c_{j}(z)b_j(z)^*}{k_j(z)},
\]
which has been considered in \cite[Equation (2.10)]{MlBeGu25}.
For example, for $\varphi^{-1}(z) = z^2 + \alpha z + \beta$ we have 
\[
G_r(z) = \bfc_r^*\left(z^2\bfI_r +  \alpha z\bfI_r - (\bfA_r-\beta\bfI_r)\right)^{-1}\bfb_r + d_r = 
\sum_{j=1}^r\frac{\phi_j}{z^2 + \alpha z + \beta -\lambda_j} +d_r,
\]
which resembles the transfer function of a second order system (see also \cite[Equation 4.2]{MlBeGu25}).
This provides an interesting connection to the $\hardy(\C_+)$ interpolatory framework introduced in \cite{MlBeGu25}.
Additionally, for a conformal map $\varphi\colon \C_+\rightarrow \A$ and $G_r$ as in \eqref{eq:Gr2} with poles in $\varphi(\lambda_p)$, and $G(z) = \sum_{j=1}^n\frac{\eta_j}{\varphi^{-1}(z)-\mu_j}$, we are able to recover interpolation conditions similar to \eqref{eq:intoptH2general}
\[
G(-\overline{\varphi(\lambda_p)}) = G_r(-\overline{\varphi(\lambda_p)}), \text{ and } G'(-\overline{\varphi(\lambda_p)}) = G_r'(-\overline{\varphi(\lambda_p)}),
\]
by using \cite[Theorem 3.3]{MlBeGu25}.

\subsubsection{Connection to the Schwarz function}\label{sec:schwarz}
In this section we discuss a connection between the interpolation points found in \cref{th:E2Dintcond} and \cref{theorem:H2Dintcond} and the concept of \textit{anti-conformal reflection}, a topic related to the Schwarz function (see \cite{Dav74}).
We start by providing a definition of Schwarz function based on \cite[Chapter 8]{Dav74}.
\begin{definition}
    Consider the arc $\Gamma$ in the complex plane as the image of the real segment $[a,b]$ under a conformal map $\psi$.
    We then locally define the Schwarz function of $\Gamma$ as
    \begin{equation}\label{eq:schwarz1}
        S(\cdot) = \overline{\psi}\circ\psi^{-1}(\cdot)    .
    \end{equation}
\end{definition}
Here, $S$ is independent of the chosen parametrization $\psi$ of $\Gamma$. 
Because our framework considers conformal maps with domain the unit disk we look at parametrizations of $\Gamma$ with respect to the unit circle. 
Let $\varphi\colon\unitdisk\rightarrow\A$ with $\Gamma=\partial\A$ a simple closed curve satisfying the following assumption based on \cite[Chapter 8]{Dav74}.
\begin{assumption}\label{assumption:schw}
    Let $\varphi$ satisfy \cref{assumption:conformalmap} and $\partial\mathbb{A}$ be an analytic simple closed curve.
    Assume $\varphi$ can be analytically continued across $\partial\unitdisk$ such that it is an analytic function in $\mathbb{X}$ with $\mathbb{X}\subseteq\{z\in\C\big| |z|<r\}$ for a fixed $r>1$.
    Further assume that $\varphi^{-1}$ can be analytically continued over $\partial\A$ in $\mathbb{Y}$ such that $\A\subset \mathbb{Y}$.   
\end{assumption}
For $\varphi$ satisfying \cref{assumption:schw}, it is possible to rewrite the Schwarz function as follows (see also \cite[Equation 8.4]{Dav74})
\[
S(\cdot) = \overline{\varphi}\left(\frac{1}{\varphi^{-1}(\cdot)}\right).
\]

In this work, we are particularly interested in the complex conjugate of the Schwarz function, known as the anti-conformal reflection $a(\cdot)=\overline{S(\cdot)}$ with respect to $\partial\A$ (see \cite[Section 1.3]{Sha92}).
If the poles $\lambda_j$ of $G_r$ are in $\mathbb{X}\backslash\unitdisk$, and \cref{assumption:schw} holds, we can rewrite them as $\lambda_j=\varphi^{-1}(\mu_j)$ with $\mu_j\in\mathbb{Y}\backslash\A$. 
It is then possible to rewrite the optimal interpolation conditions in \eqref{eq:optintE2D} and \eqref{eq:optintH2D} as follows 
        \begin{equation}\label{eq:intcondschwarz}
            \begin{aligned}
                G\circ\varphi\left(0\right) &= G_r\circ\varphi\left(0\right),\\
                G(a(\mu_j)) &= G_r(a(\mu_j)), \\
                G'(a(\mu_j))  &= G_r'(a(\mu_j)),
                \text{ for } j=1,\dots,r.
            \end{aligned}
        \end{equation}

\subsection{A TF-IRKA based algorithm}\label{sec:eTFIRKA}
In the optimization problems \eqref{eq:optE2D} and \eqref{eq:H2Dopt} the quantities used to minimize the error norm are $\phi_j$, $\lambda_j$, and $d_r$. 
We now reformulate the optimal interpolation conditions found in \cref{th:E2Dintcond} and \cref{theorem:H2Dintcond} by exploiting the structure of $G_r$ resulting in interpolation conditions for rational approximants.  
\begin{corollary}\label{corollary:optintE2Dnum}
    Let $G_r\in \etwod$ in \eqref{eq:Gr} be a local minimizer of \eqref{eq:optE2D}. We then have the following interpolation conditions
    \begin{equation}\label{eq:optintE2Dnum}    
        \begin{aligned}
            \mathcal{A}_G\left(0\right) &= -\sum_{k=1}^r \frac{\phi_k}{ \lambda_k} + d_r,\\
            \mathcal{A}_G\left(\frac{1}{\overline{\lambda_p}}\right) &= \sum_{k=1}^r \frac{\phi_k}{\frac{1}{\overline{\lambda_p}} - \lambda_k} + d_r, \\
            \frac{\mathrm{d}}{\mathrm{d}z}\left[\mathcal{A}_G(z)\right]_{z=\frac{1}{\overline{\lambda_p}}}&=
             \sum_{k=1}^r \frac{-\phi_k}{\left(\frac{1}{\overline{\lambda_p}} - \lambda_k\right)^2},
        \end{aligned}
    \end{equation}
    for $p=1,\dots,r$ and $\calA_G$ defined as in \eqref{eq:AG}.
\end{corollary}
\begin{proof}
    To prove the first two conditions in \eqref{eq:optintE2Dnum} we simply need to evaluate $G_r\circ\varphi$ at $0$ and $1/\overline{\lambda_p}$. For the third condition, first note that from the inverse function theorem, it follows that
    \begin{equation}\label{eq:neat}
        [\varphi^{-1}]'(\varphi(z)) = \frac{1}{\varphi'(z)}.
    \end{equation} 
    Hence, we obtain
    \[
        G_r'\circ\varphi(s) = \sum_{k=1}^r \frac{-\phi_k }{(s-\lambda_k)^2[\varphi'(s)]^{1/2}\varphi'(s)} - \left(\sum_{k=1}^r\frac{\phi_k}{s-\lambda_k}+d_r\right)\frac{\varphi''(s)}{2[\varphi'(s)]^{3/2}\varphi'(s)}.
    \]
    Evaluating $G_r'\circ\varphi(s)$ at $1/\overline{\lambda_p}$ results in 
    \begin{align*}
        G_r'\circ\varphi\left(\frac{1}{\overline{\lambda_p}}\right) &= \sum_{k=1}^r \frac{-\phi_k }{\left(\frac{1}{\overline{\lambda_p}}-\lambda_k\right)^2\left[\varphi'(1/\overline{\lambda_p})\right]^{1/2}\varphi'\left(1/\overline{\lambda_p}\right)} - \frac{\calA_G\left(1/\overline{\lambda_p}\right)\varphi''\left(1/\overline{\lambda_p}\right)}{2\left[\varphi'\left(1/\overline{\lambda_p}\right)\right]^{3/2}\varphi'\left(1/\overline{\lambda_p}\right)}
    \end{align*}
    where we used the second condition in \eqref{eq:optintE2Dnum}.
    We can then rewrite the evaluation as
	\begin{equation}\label{eq:eq3}
		\begin{aligned}
			G_r'\circ\varphi\left(\frac{1}{\overline{\lambda_p}}\right) 
			&= \sum_{k=1}^r \frac{-\phi_k }{\left(\frac{1}{\overline{\lambda_p}}-\lambda_k\right)^2[\varphi'\left(1/\overline{\lambda_p}\right)]^{1/2}\varphi'\left(1/\overline{\lambda_p}\right)} - \frac{G\circ\varphi\left(1/\overline{\lambda_p}\right)\frac{\mathrm{d}}{\mathrm{d}z}\left[\varphi'\left(z\right)\right]^{1/2}_{z=\frac{1}{\overline{\lambda_p}}}}{\left[\varphi'\left(1/\overline{\lambda_p}\right)\right]^{1/2}\varphi'\left(1/\overline{\lambda_p}\right)}.\\
		\end{aligned}
	\end{equation}
	Now, looking at the second condition in \eqref{eq:optintE2D} and replacing the result in \eqref{eq:eq3} results in 
	\[
	G'\circ\varphi\left(\frac{1}{\overline{\lambda_p}}\right) = \sum_{k=1}^r \frac{-\phi_k }{\left(\frac{1}{\overline{\lambda_p}}-\lambda_k\right)^2[\varphi'\left(1/\overline{\lambda_p}\right)]^{1/2}\varphi'\left(1/\overline{\lambda_p}\right)} - \frac{G\circ\varphi\left(1/\overline{\lambda_p}\right)\frac{\mathrm{d}}{\mathrm{d}z}\left[\varphi'\left(z\right)\right]^{1/2}_{z=\frac{1}{\overline{\lambda_p}}}}{\left[\varphi'\left(1/\overline{\lambda_p}\right)\right]^{1/2}\varphi'\left(1/\overline{\lambda_p}\right)}
	\]
	leading to
	\[
		G'\circ\varphi\left(\frac{1}{\overline{\lambda_p}}\right)\varphi'\left(\frac{1}{\overline{\lambda_p}}\right)\left[\varphi'\left(\frac{1}{\overline{\lambda_p}}\right)\right]^{1/2}+ G\circ\varphi\left(\frac{1}{\overline{\lambda_p}}\right)\frac{\mathrm{d}}{\mathrm{d}z}\left[\varphi'\left(z\right)\right]^{1/2}_{z=\frac{1}{\overline{\lambda_p}}}=\sum_{k=1}^r \frac{-\phi_k }{(\frac{1}{\overline{\lambda_p}}-\lambda_k)^2} 
	\]
	proving the last equality in \eqref{eq:optintE2Dnum}.
\end{proof}

\begin{corollary}\label{corollary:optintH2Dnum}
    Let $G_r\in \durenhardy$ in \eqref{eq:Gr2} be a local minimizer of \eqref{eq:H2Dopt}. We then have the following interpolation conditions
    \begin{equation}\label{eq:optintH2Dnum}    
        \begin{aligned}
            G\circ\varphi(0) &= -\sum_{j=1}^r \dfrac{\phi_j}{\lambda_j} + d_r,\\
            G\circ\varphi\left(\frac{1}{\overline{\lambda_p}}\right) &= \sum_{k=1}^r \frac{\phi_k}{\frac{1}{\overline{\lambda_p}} - \lambda_k}+d_r,\\
            G'\circ\varphi\left(\frac{1}{\overline{\lambda_p}}\right)\varphi'\left(\frac{1}{\overline{\lambda_p}}\right) &= \sum_{k=1}^r \frac{-\phi_k}{\left(\frac{1}{\overline{\lambda_p}} - \lambda_k\right)^2},
        \end{aligned}   
    \end{equation}
    for $p=1,\dots,r$.
\end{corollary}
\begin{proof}
The first two conditions can be proved in the same way as in \cref{corollary:optintE2Dnum}. 
We have that the derivative of $G_r$ is given by
\[
G_r'(z) = \sum_{j=1}^r \frac{-\phi_j}{(\varphi^{-1}(z)-\lambda_j)^2}[\varphi^{-1}]'(z)
\]
which, when evaluated at $\varphi(1/\overline{\lambda_p})$ leads to 
\begin{equation}\label{eq:Grpphi}
    G_r'\circ\varphi\left(\frac{1}{\overline{\lambda_p}}\right) = \sum_{j=1}^r \frac{-\phi_j}{\left(\frac{1}{\overline{\lambda_p}}-\lambda_j\right)^2}[\varphi^{-1}]'\left(\varphi\left(\frac{1}{\overline{\lambda_p}}\right)\right) = \sum_{j=1}^r \frac{-\phi_j}{\left(\frac{1}{\overline{\lambda_p}}-\lambda_j\right)^2}\left(\varphi'\left(\frac{1}{\overline{\lambda_p}}\right)\right)^{-1}.
\end{equation}
When \eqref{eq:Grpphi} is replaced in the last equality of \eqref{eq:optintH2D} we then get the last condition in \eqref{eq:optintH2Dnum}.
\end{proof}

\begin{figure}[hbt]
    \centering
    \begin{tikzpicture}[
            node distance=1cm and 3cm,
            square/.style={draw, minimum size=2cm, align=center},
            every node/.style={font=\small}
        ]
        \node[square] (A) {$G$};
        \node[square, right=of A] (B) {$G\circ\varphi$\\$\calA_{G}$};
        \node[square, below=of A] (C) {$G_r$};
        \node[square, below=of B] (D) {$\displaystyle \sum_{j=1}^r\dfrac{\phi_j}{\cdot-\lambda_j}$};
        \draw[->, thick] (A) -- (B) node[midway, above left] {$\varphi$};
        \draw[<-, thick] (C) -- (D) node[midway, below right] {$\varphi^{-1}$};
        \draw[->, thick] (A) -- (C) node[midway, left] {\cref{alg:eTFIRKA}};
        \draw[->, thick] (B) -- (D) node[midway, right] {TF-IRKA};

        \draw[dashed] ($($(A)!0.5!(B)$) + (0,2)$) -- ($($(C)!0.5!(D)$) + (0,-1.5)$);

        \node[above=1.3cm of $(A)!0.1!(B)$] {\large $\etwod$ or $\durenhardy$};
        \node[above=1.3cm of $(A)!0.8!(B)$] {\large $\hardy$};

    \end{tikzpicture}
    \caption{Illustration depicting the steps of \cref{alg:eTFIRKA}.}
    \label{fig:illustrationalg1}
\end{figure}
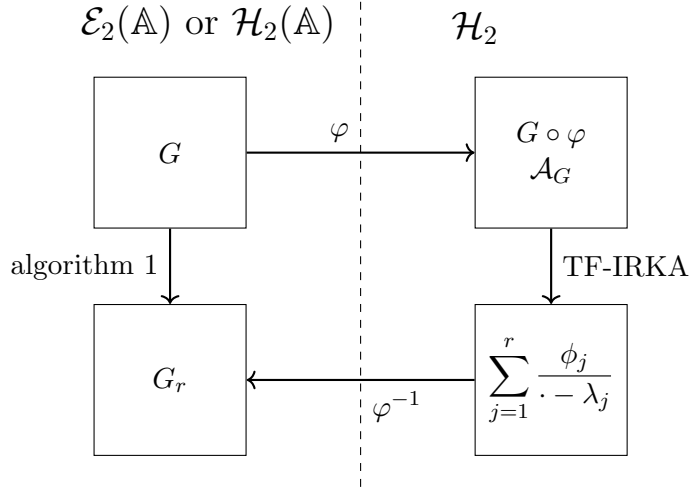

Now that the approximant has a rational structure, we can use a variation of TF-IRKA from \cite{BeGu12} for discrete-time systems to approximate $\mathcal{A}_G$ for the $\etwod$ framework and $G\circ\varphi$ for the $\durenhardy$ framework.
While in \cite{BeGu12} the $\hardy(\C_+)$ optimal interpolation conditions are used as an update for the interpolation points, we use the same algorithm but with $1/\overline{\lambda}$ as interpolation points to satisfy \cref{corollary:optintE2Dnum} and \cref{corollary:optintH2Dnum}.
Once $\phi_j$, $\lambda_j$, and $d_r$ are computed, we reconstruct $G_r$ based on \eqref{eq:Gr} or \eqref{eq:Gr2}.
A pseudocode of the algorithm is presented in \cref{alg:eTFIRKA}.
Additionally, we provide an illustration of \cref{alg:eTFIRKA} in \cref{fig:illustrationalg1} which takes inspiration from \cite{ClCrAl96}.
In \cref{alg:eTFIRKA} we indicate the function to approximate as $F$, which in our case is either $\mathcal{A}_G$ or $G\circ\varphi$ depending on the space of functions used.
To account for the term $d_r$ we follow \cite[Theorem 4.2.3]{AnBeGu20}
and, given the points $\{\sigma_i\}_{i=1}^r$, set 
\begin{equation}\label{eq:coeffGr}
	\bfE_r = -\mathbf{L}, \; \bfA_r = d_r - \mathbf{M}, \; \bfb_r = \mathbf{Z}^\top - d_r, \; \bfc_r = \mathbf{Y} - d_r,
\end{equation}
where $\mathbf{Z} = \mathbf{Y} = \begin{bmatrix} F(\sigma_1), \dots, F(\sigma_r) \end{bmatrix}$, and 
\begin{equation*}
	\begin{aligned}
		\mathbf{L}_{ij} &= \begin{cases}
			&\dfrac{F(\sigma_i) - F(\sigma_j)}{\sigma_i-\sigma_j}, \text{ if } i\neq j,\\
			&F'(\sigma_i), \text{ if } i=j,
		\end{cases}, \quad
		\mathbf{M}_{ij} = \begin{cases}
			&\dfrac{\sigma_iF(\sigma_i) - \sigma_j F(\sigma_j)}{\sigma_i-\sigma_j}, \text{ if } i\neq j, \\
			&F(\sigma_i) + \sigma_iF'(\sigma_i), \text{ if } i=j.
		\end{cases}
	\end{aligned}   
\end{equation*}
\begin{algorithm}
    \begin{algorithmic}[1]
    \Require{the function to approximate $F$, conformal map $\varphi$, reduced order $r$, initial interpolation points $\sigma^{(0)}$, variable $\texttt{framework}$ indicating either $\etwod$ or $\durenhardy$}
    \While{$\|\sigma^{(i+1)}-\sigma^{(i)}\|/\|\sigma^{(i)}\|>\texttt{tol}$}
        \State $d_r = F(0) + \bfc_r^\top \bfA_r^{-1}\bfb_r$
        \State Construct $\bfE_r$, $\bfA_r$, $\bfb_r$, and $\bfc_r$ based on \eqref{eq:coeffGr}
        \State Compute the eigenvalues $\bfA_r \mathbf{x}_j = \lambda_j\bfE_r\mathbf{x}_j$ for $j=1,\dots,r$
        \State $\sigma_{j}^{(i+1)}=1/\overline{\lambda_j}$ for $j=1,\dots,r$
    \EndWhile
    \If{$\texttt{framework}==\etwod$}

        \Return $G_r(z) = \left(\bfc_r^\top(\varphi^{-1}(z)\bfE_r-\bfA_r)^{-1}\bfb_r + d_r\right)[\varphi'(\varphi^{-1}(z))]^{-1/2} $
    
    \ElsIf{$\texttt{framework}==\durenhardy$}
        
        \Return $G_r(z) = \bfc_r^\top(\varphi^{-1}(z)\bfE_r-\bfA_r)^{-1}\bfb_r + d_r$
    \EndIf
    
    \caption{TF-IRKA on general domains} \label{alg:eTFIRKA}
    \end{algorithmic}
\end{algorithm}
\begin{remark}\label{remark:difficulty}
    The function $F$ to approximate on $\unitdisk$ is either $\mathcal{A}_G$ or $G\circ\varphi$, which can result in a more computationally challenging task than the approximation of $G$ itself. 
    The choice of $\varphi$ might dictate the difficulty of the approximation.
\end{remark}

\subsection{Model order reduction of discrete-time time invariant delay systems}\label{sec:mor_disc}
Let us consider a conformal map of the form $\varphi(z)=p(z)/q(z)$ where $p(z) = \sum_{j=0}^{n_p}\alpha_jz^j$ and $q(z) = \sum_{i=0}^{n_q}\beta_iz^i$ are polynomials of order $n_p\geq 1$ and $n_q\geq 1$, respectively, with $\{\alpha_j\}_{j=0}^{n_p}\in\C$ and $\{\beta_i\}_{i=0}^{n_q}\in\C$. 
Let $G(\cdot) = \bfc^*(\cdot\bfI-\bfA)^{-1}\bfb\in\durenhardy$ be the transfer function of an LTI system. 
We then look at the composition with the conformal map $\varphi$
\begin{equation*}
    G\circ\varphi(\cdot)  = \bfc^* \left(\varphi(\cdot)\bfI - \bfA\right)^{-1} \bfb = \bfc^*  \left(p(\cdot)\bfI - q(\cdot)\bfA\right)^{-1} q(\cdot)\bfb \in\hardy. 
\end{equation*}
Due to $G\circ\varphi$ being analytic in $\unitdisk$, we make the connection to dynamical systems by using the transformation $z\rightarrow z^{-1}$ (see \cite[Section 3.2.3]{ChuChe97}) leading to
\begin{equation}\label{eq:tftds}
    \begin{aligned}
        G\circ\varphi(z^{-1}) 
        &= \bfc^*  \left(\sum_{j=0}^{n_p}\alpha_jz^{-j}\bfI - \sum_{i=0}^{n_q}\beta_iz^{-i}\bfA\right)^{-1}\left(\sum_{i=0}^{n_q}\beta_iz^{-i}\right)\bfb
    \end{aligned}
\end{equation}
which is the transfer function of the following discrete-time delay system 
\begin{equation}\label{eq:tds1}
    \begin{cases}
        \sum_{j=0}^{n_p}\alpha_j x_{k-j} - \sum_{i=0}^{n_q}\beta_i\bfA x_{k-i} = \bfb \sum_{i=0}^{n_q}\beta_i u_{k-j},\\
        y_k = \bfc^* x_k,
    \end{cases}
\end{equation}
for $k\in\mathbb{N}$ and initial conditions $x(\ell)=0$ for $\ell=-\max(n_p,n_q),\dots,0$, cp.~\cite[Chapter 6.1.2]{Fri2014}.

Similar to \cite{GuAnBe08}, there exists a bound on the $\ell_\infty$ error between the output of \eqref{eq:tds1} and the reduced model $G_r\circ\varphi$. 
More in detail, let $\widehat{y}_k$ be the output at step $k$ of the LTI system 
\begin{equation}\label{eq:ssGr_disc}
    \begin{cases}
        \bfA_r \widehat{x}_{k} = \bfE_r \widehat{x}_{k-1} - \bfb_r u_{k}, \quad \widehat{x}_{\ell} = 0 \text{ for } \ell=\{-1,0\},\\
        \widehat{y}_k = \bfc_r^* \widehat{x}_k + d_ru_k,
    \end{cases}
\end{equation}
corresponding to the inverse $\mathcal{Z}$-transform of 
\begin{equation*}
    G_r\circ\varphi(z^{-1}) \longrightarrow
    \begin{cases}
        \left(z^{-1}\bfE_r - \bfA_r\right)\widehat{X}(z) = \bfb_rU(z),\\
        \widehat{Y}(z) = \bfc_r^* \widehat{X}(z) + d_r U(z).
    \end{cases}
\end{equation*}
where $\bfA_r\in\C^{r\times r}$, $\bfb_r,\bfc_r\in\C^r$, $d_r\in\C$.
Denote by $Y$ and $U$ the inverse $\mathcal{Z}$-transform of $y_k$ and $u_k$ in \eqref{eq:tds1}.
We then get
\begin{equation}\label{eq:linftybound}
    \begin{aligned}
    \|y-\widehat{y}_k\|_{\ell_\infty} &= \max_{k>0}|y_k-\widehat{y}_k| = \max_{k>0}\left|\frac{1}{2\pi\iunit}\oint_{\partial\unitdisk}\left(Y(z)-\widehat{Y}(z)\right)z^{k-1}\mathrm{d}z \right| \\
    &= \max_{k>0}\left|\frac{1}{2\pi}\int_{0}^{2\pi}\left(Y(e^{\iunit\theta})-\widehat{Y}(e^{\iunit\theta})\right)e^{\iunit\theta k}\mathrm{d}\theta \right| \leq \frac{1}{2\pi}\int_{0}^{2\pi}\left|Y(e^{\iunit\theta})-\widehat{Y}(e^{\iunit\theta})\right|\mathrm{d}\theta\\
    &= \frac{1}{2\pi}\int_{0}^{2\pi}\left|G\circ\varphi(e^{-\iunit\theta})-G_r\circ\varphi(e^{-\iunit\theta})\right|\left|U(e^{\iunit\theta})\right|\mathrm{d}\theta\\
    &\leq\left(\frac{1}{2\pi}\int_{0}^{2\pi}\left|G\circ\varphi(e^{-\iunit\theta})-G_r\circ\varphi(e^{-\iunit\theta})\right|^2\mathrm{d}\theta\right)^{1/2}\left(\frac{1}{2\pi}\int_{0}^{2\pi}\left|U(e^{\iunit\theta})\right|^2\mathrm{d}\theta\right)^{1/2}\\
    &= \|G\circ\varphi-G_r\circ\varphi\|_{\hardy}\|U\circ e^{\iunit\cdot}\|_{\mathcal{L}_2(0,2\pi)} = \|G-G_r\|_{\durenhardy}\|u_k\|_{\ell_2},
    \end{aligned}
\end{equation}
where the last equality results from Parseval's theorem.

This theoretical framework also applies to discretized LTI systems. 
Here, the choice of $\varphi$ coincides with the chosen discretization technique. 
Below we provide some examples and the resulting interpolation conditions for asymptotically stable systems with transfer function $G(z) = \bfc^\top(z\bfI-\bfA)^{-1}\bfb$. 
Connections between $\hardy(\C_+)$ model order reduction and discretization methods have also been discussed in \cite[Sections 6.2, 6.3, and 6.4]{Lat16}, however, without any link to the optimal interpolation conditions in \cref{theorem:H2Dintcond}.
Further details on the connections to \cite{Lat16} are given below.

\subsubsection{Implicit Euler}\label{sec:impliciteuler}
\begin{figure}[hbt]
    \centering
    \input{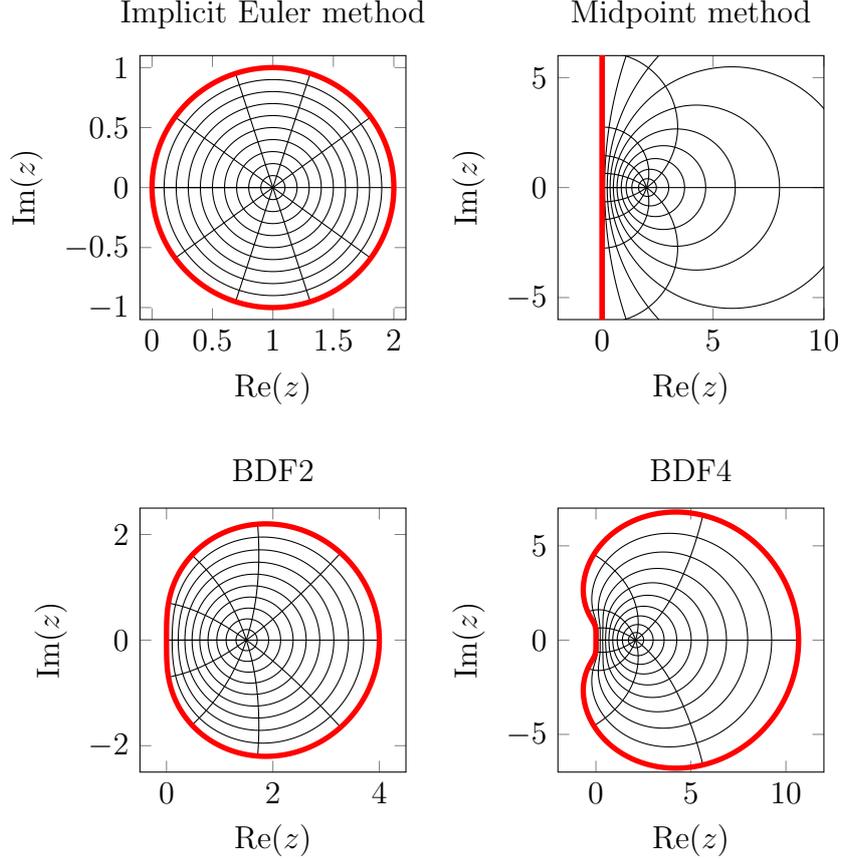}
    \caption{The image $\A$ of $\varphi$ for different discretization techniques with time step $h=1$.
    The red contour indicates $\partial\A$.}
    \label{fig:mappings}
\end{figure}

Applying the implicit Euler method to an LTI system characterized by $G$ yields
\begin{equation}\label{eq:impleuler}
        \begin{cases}
            x_{k+1} = x_k + h\bfA x_{k+1} + h\bfb u_{k+1}, \quad x_0=0,\\
            y_k = \bfc^\top x_k,
        \end{cases}
\end{equation}
where $x_k = x(kh)$ and $h$ denotes the chosen time step.
Applying the $\mathcal{Z}$-transform yields
\[
\begin{cases}
    (z\bfI - zh\bfA - \bfI)X(z) = zh\bfb U(z),\\
    Y(z) = \bfc^\top X(z),    
\end{cases}
\]
with transfer function 
\[
H(z) = \bfc^\top \left((z-1)\bfI - zh\bfA\right)^{-1}zh\bfb = \bfc^\top \left(\frac{1-z^{-1}}{h}\bfI - \bfA\right)^{-1}\bfb.
\]
Consequently, if we choose 
\begin{equation}\label{eq:confmapimpeul}
\varphi(z) = \frac{1-z}{h}
\end{equation}
it holds that $G\circ\varphi(z^{-1})=H(z)$.
Let us emphasize that \eqref{eq:confmapimpeul} conformally maps the unit disk into a disk of center $1/h$ and radius $1/h$ (see \cref{fig:mappings} for $h=1$), i.e., $\A=\left\{z\in\C\;| \;|z-1/h|<1/h\right\}$. 
For a time step of $h=1$, the interpolation conditions in \cref{theorem:H2Dintcond} then become
\begin{equation}\label{eq:impliciteulerintcond}
    \begin{aligned}
        G\left(1\right) = G_r\left(1\right), \; G\left(1-\frac{1}{\overline{\lambda_p}}\right) = G_r\left(1-\frac{1}{\overline{\lambda_p}}\right), \text{and }  G'\left(1-\frac{1}{\overline{\lambda_p}}\right) = G_r'\left(1-\frac{1}{\overline{\lambda_p}}\right), 
    \end{aligned}
\end{equation}
for $p=1,\dots,r$.
Due to the conformal map in \eqref{eq:confmapimpeul} being bijective and conformal in $\C$, it satisfies \cref{assumption:schw}.  
We can thus rewrite the interpolation conditions through the Schwarz function introduced in \cref{sec:schwarz}.
Let $\mu_p = \varphi^{-1}(\lambda_p)$, using \eqref{eq:intcondschwarz} with $h=1$ we then get 
\begin{equation}\label{eq:impliciteulerintcond_schwarz}
    \begin{aligned}
        G\left(1\right) = G_r\left(1\right), \; G\left(\frac{\overline{\mu_p}}{\overline{\mu_p}-1}\right) = G_r\left(\frac{\overline{\mu_p}}{\overline{\mu_p}-1}\right), \; \text{and }  G'\left(\frac{\overline{\mu_p}}{\overline{\mu_p}-1}\right) = G_r'\left(\frac{\overline{\mu_p}}{\overline{\mu_p}-1}\right),
    \end{aligned}
\end{equation}
for $p=1,\dots,r$.
In \cref{fig:anticonformal} we show the result of applying the anti-conformal reflection to points in $\A$.
It is easy to see that, due to the affine nature of $\varphi$, $G_r$ will still be a rational function with poles in $1-h\lambda$
\begin{equation}\label{eq:GrimplicitEuler}
G_r(z) = \bfc_r^*\left(\varphi^{-1}(z)\bfI_r - \mathbf{\Lambda}_r\right)^{-1}\bfb_r +d_r = -\bfc_r^*\left(z\bfI_r - \left(\bfI_r-h\mathbf{\Lambda}_r\right)\right)^{-1}\frac{1}{h}\bfb_r + d_r.
\end{equation}

\begin{remark}
In \cite[Theorem 6.20]{Lat16} it is shown that a ROM satisfying the $\hardy(\C_+)$ optimality conditions in \eqref{eq:intoptH2Cplus} does not retain the $\hardy(\bar{\unitdisk}^{\mathsf{c}})$ optimality conditions in \eqref{eq:intoptH2} when the system is discretized using the implicit Euler method. 
In our framework, $G_r$ being $\hardy(\C_+)$ optimal would require $\A=\C_+$. 
However, applying the implicit Euler method instead results in an affine conformal map with image $\A=\left\{z\in\C\;| \;|z-1/h|<1/h\right\}$, a disk in the complex plane. The result provided in \cite[Theorem 6.20]{Lat16} thus translates to stating that $\varphi(z)=(1-z)/h$, appearing in the implicit Euler discretization, does not map the unit disk $\unitdisk$ into the right half plane $\C_+$.
As a remedy, our framework assumes $G_r$ to be $\durenhardy$ optimal and satisfying \eqref{eq:impliciteulerintcond_schwarz} rather than assuming $\hardy(\C_+)$-optimality from the start.
\end{remark}

\subsubsection{Midpoint method}\label{sec:trapezoidal}
We start with the discretized version of the LTI system through the midpoint method
\begin{equation}
    \begin{cases}
        x_{k+1} = x_k + \frac{h}{2}\bfA\left(x_k+x_{k+1}\right) + \frac{h}{2}\bfb \left(u_k+u_{k+1}\right), \quad x(0)=0\\
        y_k = \bfc^\top x_k.
    \end{cases}
\end{equation}
Applying the $\mathcal{Z}$-transform leads to the transfer function
\[
H(z) = \bfc^\top \left((z-1)\bfI - (z+1)\frac{h}{2}\bfA \right)^{-1} (z+1)\frac{h}{2}\bfb = \bfc^\top \left(\frac{2}{h}\frac{z-1}{z+1}\bfI - \bfA \right)^{-1}\bfb = G\circ\varphi(z^{-1}).
\]
Here, the conformal map coincides with the M\"obius transform
\begin{equation}\label{eq:tustin}
    \varphi(z) = \frac{2}{h}\frac{1-z}{1+z},
\end{equation}
which maps the unit disk to the right half plane (see \cref{fig:mappings}). 
The optimality conditions in \eqref{theorem:H2Dintcond} translate to
\begin{equation}\label{eq:trapezoidalintcond}
    \begin{aligned}
        G\left(2/h\right) = G_r\left(2/h\right), \; G\left(\frac{2}{h}\frac{\overline{\lambda}_p-1}{\overline{\lambda}_p+1}\right) &= G_r\left(\frac{2}{h}\frac{\overline{\lambda}_p-1}{\overline{\lambda}_p+1}\right),\\ \text{and }  G'\left(\frac{2}{h}\frac{\overline{\lambda}_p-1}{\overline{\lambda}_p+1}\right) &= G_r'\left(\frac{2}{h}\frac{\overline{\lambda}_p-1}{\overline{\lambda}_p+1}\right), 
    \end{aligned}
\end{equation}
for $p=1,\dots,r$.
Also in this case we can rewrite \eqref{eq:trapezoidalintcond} using the anti-conformal reflection introduced in \cref{sec:schwarz}.
Let $\lambda_p=\varphi^{-1}(\mu_p) = \frac{1-h\mu_p/2}{1+h\mu_p/2}$ with $\mu_p\in\C_-$, we then get 
\begin{equation}\label{eq:trapezoidalintcond_shcwarz}
    \begin{aligned}
        G\left(2/h\right) = G_r\left(2/h\right), \; G\left(-{\overline{\mu_p}}\right) = G_r\left(-{\overline{\mu_p}}\right), \; \text{and }  G'\left(-{\overline{\mu_p}}\right) = G_r'\left(-{\overline{\mu_p}}\right),
    \end{aligned}
\end{equation}
for  $p=1,\dots,r$.
It is possible to see that, for $d_r=0$, \eqref{eq:trapezoidalintcond_shcwarz} resembles the $\hardy(\C_+)$ optimality conditions in \eqref{eq:intoptH2Cplus}. 
A similar result can be found in \cite[Theorem 6.21]{Lat16} where it was proven that the resulting $\hardy(\C_+)$ optimal ROM preserves the $\hardy(\bar{\unitdisk}^{\mathsf{c}})$ optimality conditions when discretized with the midpoint method.  
This is due to the equivalent conformal map being the M\"obius transform \eqref{eq:tustin} that maps the unit disk into $\A=\C_+$, resulting in the preservation of the interpolation conditions between the two Hardy spaces.

Once more, $G_r$ has a rational structure due to $\varphi^{-1}$ being rational itself.
More in detail, we have
\begin{equation}\label{eq:rom_tustin}
    G_r(z) = \bfc_r^*\left(\varphi^{-1}(z)\bfI_r - \mathbf{\Lambda}_r\right)^{-1} \bfb_r + d_r  = \widehat{\bfc}_r^*\left(z\bfI_r - \widehat{\mathbf{\Lambda}}_r\right)^{-1} \widehat{\bfb}_r + \widehat{d}_r,
\end{equation}
where $\widehat{\mathbf{\Lambda}}_r = \frac{2}{h}\left(\bfI_r + \mathbf{\Lambda}_r\right)^{-1}\left(\bfI_r - \mathbf{\Lambda}_r\right)$, $\widehat{\bfc}_r^* = \bfc_r^*\left(\widehat{\mathbf{\Lambda}}_r + \frac{2}{h}\bfI_r\right)$, $\widehat{\bfb}_r = -\left(\widehat{\mathbf{\Lambda}}_r + \frac{2}{h}\bfI_r\right)^{-1}\bfb_r$, and $\widehat{d}_r = \bfc_r^\top\widehat{\bfb}_r+d_r$. Note that the poles of \eqref{eq:rom_tustin} are $\varphi(\lambda_j)\in\C_-$. 
Considering $G$ being asymptotically stable, an $\durenhardy$ optimal ROM satisfying \eqref{eq:trapezoidalintcond} will preserve the stability of the original FOM.

\subsubsection{BDF methods}\label{sec:bdf}
\begin{figure}[hbt]
    \centering
    \input{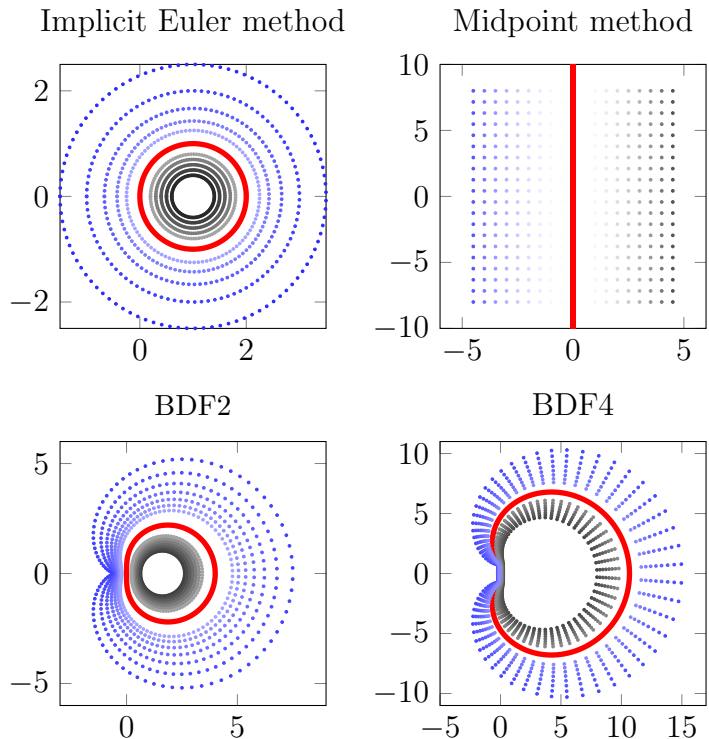}
    \caption{Depiction of the anti-conformal reflection with respect to $\partial\A$ for different discretization techniques. 
    The grey points are mapped into the blue ones.}
    \label{fig:anticonformal}
\end{figure}

In the following we show that also the BDF2 and BDF4 methods fit into our framework. 
Similar derivations can be applied to all other BDF methods. 
\paragraph*{BDF2}
Let us consider the discretized LTI system with the BDF2 method
\begin{equation}\label{eq:ssbdf2}
    \begin{cases}
        x_{k+2} - \frac{4}{3}x_{k+1} + \frac{1}{3}x_{k} = \frac{2}{3}h\left(\bfA x_{k+2} + \bfb u_{k+2}\right),\\
        y_k = \bfc^\top x_k,
    \end{cases}
\end{equation}
with $x(0)=0$. 
Applying the $\mathcal{Z}$-transform leads to 
\begin{equation*}
    \begin{cases}
        \left((\bfI - \frac{2h}{3}\bfA)z^2  - \frac{4}{3}z + \frac{1}{3}\right)X(z) = \frac{2}{3}h\bfb z^2 U(z) ,\\
        Y(z) = \bfc^\top X(z),
    \end{cases}
\end{equation*}
with transfer function
\begin{equation*}
H(z) = \bfc^\top \left(\frac{1}{h}\left(\frac{z^{-2}}{2} - 2z^{-1} + \frac{3}{2}\right)\bfI - \bfA\right)^{-1} \bfb = G\circ\varphi(z^{-1}),
\end{equation*}
resulting in
\begin{equation}\label{eq:bdf2confmap}
    \varphi(z) = \frac{1}{h}\left(\frac{z^{2}}{2} - 2z + \frac{3}{2}\right) \text{ and } \varphi^{-1}(z) = 2-\sqrt{2hz+1}.
\end{equation}
The image $\A$ of $\varphi$ is depicted in \cref{fig:mappings}.
From $\varphi^{-1}$ in \eqref{eq:bdf2confmap} it is clear that $G_r$ does not have a rational structure.
While the optimality conditions in \eqref{eq:optintH2D} are straightforward to derive with \eqref{eq:bdf2confmap}, we look at the Schwarz function variation in \eqref{eq:intcondschwarz}.
In this case we consider points $\lambda_p$ and $\mu_p$ in neighborhoods of $\partial\unitdisk$ and $\partial\A$, respectively. 
Let $\lambda_p=\varphi^{-1}(\mu_p)$, we then get \eqref{eq:intcondschwarz} with 
\begin{equation*}\label{eq:bdf2intcond_shcwarz}
    a(z) = \frac{1}{2h}\overline{\left(\frac{6hz - 8\sqrt{2hz+1} + 8}{2hz - 4\sqrt{2hz + 1} + 5}\right)}.
\end{equation*}
The effect of the anti-conformal function \eqref{eq:bdf2intcond_shcwarz} to points inside $\partial\A$ is illustrated in \cref{fig:anticonformal}.

\paragraph*{BDF4}
Let us consider the discretized LTI system with the BDF4 method
\begin{equation}\label{eq:ssdiscbdf4}
    \begin{cases}
        x_{k+4} - \frac{48}{25}x_{k+3} + \frac{36}{25}x_{k+2} - \frac{16}{25}x_{k+1} + \frac{3}{25}x_k = \frac{12}{25}h\left(\bfA x_{k+4} + \bfb u_{k+4}\right),\\
        y_k = \bfc^\top x_k,
    \end{cases}
\end{equation}
with $x(0)=0$. 
Applying the $\mathcal{Z}$-transform leads to 
\begin{equation*}
    \begin{cases}
        \left((\bfI - \frac{12}{25}h\bfA)z^4  - \frac{48}{25}z^3 + \frac{36}{25}z^2 - \frac{16}{25}z + \frac{3}{25}\right)X(z) = \frac{12}{25}h\bfb z^4 U(z),\\
        Y(z) = \bfc^\top X(z),
    \end{cases}
\end{equation*}
with transfer function
\begin{equation*}
H(z) = \bfc^\top \left(\frac{1}{h}\left(\frac{z^{-4}}{4} - \frac{4}{3}z^{-3} + 3z^{-2} - 4z^{-1} + \frac{25}{12}\right)\bfI - \bfA\right)^{-1} \bfb = G\circ\varphi(z^{-1}),
\end{equation*}
where 
\begin{equation}\label{eq:confmapbdf4}
\varphi(z) = \frac{1}{h}\left(\frac{z^{4}}{4} - \frac{4}{3}z^{3} + 3z^{2} - 4z + \frac{25}{12}\right),
\end{equation}
is the conformal map.
We do not show the inverse of $\varphi$ due to its convoluted formulation.
Nevertheless, in \cref{fig:anticonformal} we show the effect of the resulting anti-conformal map $a$ approximated with the adaptive Antoulas Anderson algorithm (for more details see \cite{Tre25}).
Similarly to BDF2, also here $G_r$ does not have a rational structure.

\section{Numerics}\label{sec:numerics}
In this section, we test \cref{alg:eTFIRKA} on three numerical examples. 
In \cref{sec:numCplus} we use the $\etwod$ framework to approximate the transfer function of the one-dimensional heat equation and compare the results with TF-IRKA from \cite{BeGu12}.
In \cref{sec:numbdf2} we look at the clamped beam from the SLICOT library discretized with the BDF2 method. 
Here, we use the $\durenhardy$ framework and test the performances of \cref{alg:eTFIRKA} against the adaptive Antoulas-Anderson algorithm. 
Lastly, in \cref{sec:bdf4num} we reduce a two-dimensional convection-diffusion equation spatially discretized with finite differences and discretized in time with BDF4. 

All the experiments were carried out in MATLAB 2023b on an Apple Macbook Pro with an Apple M2 Pro CPU and 16GB of RAM. 
In all experiments we initialize \cref{alg:eTFIRKA} with the following shifts
$$\sigma^{(0)} = 0.1*\texttt{randn}(r,1) + 0.1\iunit*\texttt{randn}(r,1); \; \sigma^{(0)} = \frac{\sigma^{(0)}}{|\sigma^{(0)}|}*\texttt{rand}(r,1),$$ 
where we used the MATLAB commands $\texttt{randn}$ and $\texttt{rand}$, and both division and multiplications are element-wise.
Additionally, we set $\texttt{tol}=10^{-6}$ for \cref{alg:eTFIRKA}.
Both the $\etwod$ and $\durenhardy$ norms are computed with the \texttt{integral} command of MATLAB where the relative and absolute tolerances are set to $10^{-8}$ and $10^{-12}$, respectively. 
For AAA we use the standard tolerance of $10^{-13}$. 
Lastly, the system matrices in \cref{sec:bdf4num} are constructed with the functions $\texttt{fdm\_2d\_matrix}$ and $\texttt{fdm\_2d\_vector}$ from \cite{Pen99}.
The code can be found in \url{https://github.com/aaborghi/dd-H2-simplyconnected.git}.

\subsection{Half plane region}\label{sec:numCplus}
We use the transfer function of the one-dimensional heat equation describing the temperature distribution on a semi-infinite rod given in \cite{BeGu12}.
The transfer function of this dynamical system results in
\[
    G(z) = e^{-\sqrt{z}}.
\]
To conform our numerical experiment with \cite{BeGu12} we choose $\A=\C_+$ and use the M\"obius transform
\begin{equation}\label{eq:moebius}
    \varphi(z) = \frac{1+z}{1-z}, \quad \varphi'(z) = \frac{2}{(1-z)^2}.
\end{equation}

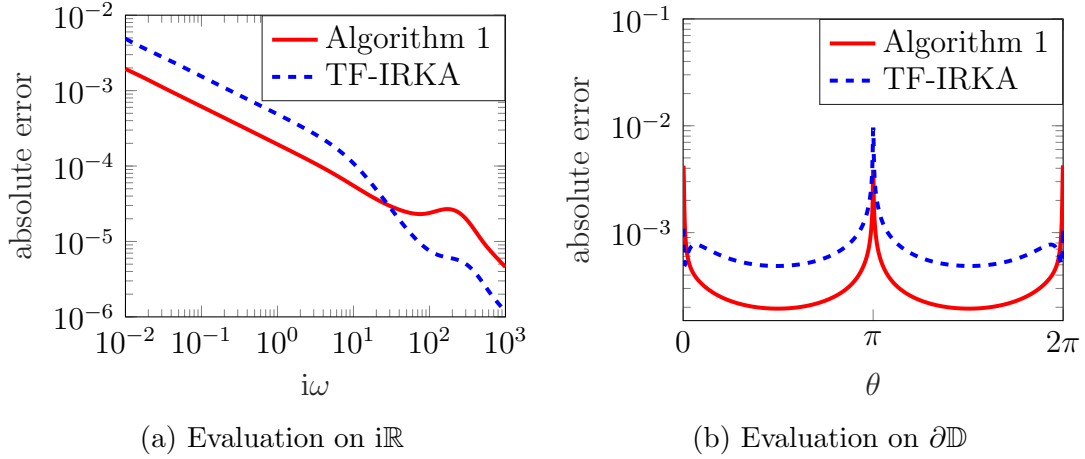
\begin{figure}[hbt]
    \centering 
    \begin{subfigure}[b]{0.45\textwidth}
        \input{figures/iR_E2Cplus.tex}
        \caption{Evaluation on $\iunit\R$}
        \label{fig:E2Cplus_a}
    \end{subfigure}
    \begin{subfigure}[b]{0.45\textwidth}
        \input{figures/partialD_E2Cplus.tex}
        \caption{Evaluation on $\partial\unitdisk$}
        \label{fig:E2Cplus_b}
    \end{subfigure}
    \caption{The absolute error between the transfer function of the heat equation and the approximant computed by TF-IRKA and \cref{alg:eTFIRKA} for $r=10$.}
    \label{fig:E2Cplus}
\end{figure}

\begin{figure}[hbt]
    \centering 
    \begin{subfigure}[b]{0.45\textwidth}
        \includegraphics[width=1.1\textwidth]{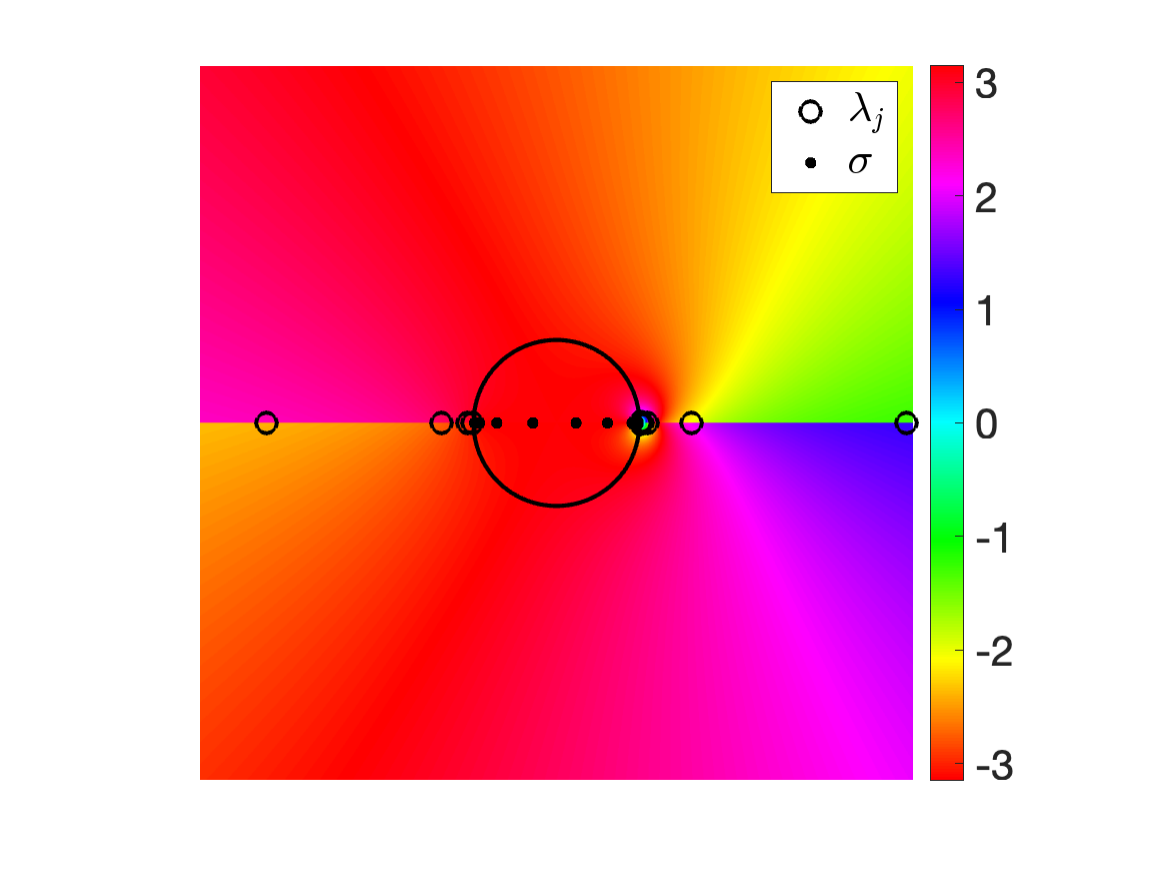}
        \caption{Phase portrait}
        \label{fig:halfplane_a}
    \end{subfigure}
    \begin{subfigure}[b]{0.45\textwidth}
        \input{figures/E2Cplus.tex}
        \caption{$\etwod$ relative error}
        \label{fig:halfplane_b}
    \end{subfigure}
    \caption{(Left) Phase plot of $\calA_G$ with the poles $\lambda_j$ of the rational approximant and the final shifts $\sigma$ computed by \cref{alg:eTFIRKA} for $r=10$.
    (Right) $\etwod$ relative error norm between $\calA_G$ and the approximants computed with \cref{alg:eTFIRKA} and TF-IRKA for different values of $r$.}
    \label{fig:halfplane}
\end{figure}

In \cref{fig:E2Cplus} we show a comparative error plot of the resulting ROM for $r=10$ computed with TF-IRKA developed in \cite{BeGu12} and \cref{alg:eTFIRKA} with the $\etwod$ framework.
In \cref{fig:E2Cplus_a} the absolute error on the imaginary axis of both reduced models is displayed. 
More in detail, we define the absolute error to be $|G(\iunit\omega) - F_r(\iunit\omega)|$ for $\omega\in\R$ where $F_r(\cdot)=\bfc_r^*\left(\varphi^{-1}(\cdot)\bfE_r - \bfA_r\right)^{-1}\bfb_r$ when \cref{alg:eTFIRKA} is applied, and $F_r(\cdot) = \tilde{\bfc}_r^*(\cdot\tilde{\bfE}_r - \tilde{\bfA}_r)^{-1}\tilde{\bfb}_r$ when TF-IRKA is used. 
While \cref{alg:eTFIRKA} provides a lower error until $\omega\approx 12$, TF-IRKA provides better results at higher frequencies.
A similar conclusion can be made in \cref{fig:E2Cplus_b} where the absolute error $|\mathcal{A}_G(e^{\iunit\theta})-\mathcal{A}_{F_r}(e^{\iunit\theta})|$ over the unit circle is illustrated.
Here, \cref{alg:eTFIRKA} provides a lower error in most of the unit circle, with a deteriorating performance towards $0$ and $2\pi$.
This is due to $\theta\rightarrow 0$ and $\theta\rightarrow 2\pi$ coinciding with $\omega\rightarrow\pm \infty$. 
The resulting pole placement of \cref{alg:eTFIRKA} and the corresponding shifts $\sigma$ are shown in \cref{fig:halfplane_a} with the phase plot of $\calA_G$. 
The discontinuity present in the phase plot is a result of $G$ having a branch cut on the negative real axis. 
Under the M\"obious transform the negative real line is mapped into the line from $1$ to $-1$ that passes through infinity. 
It is possible to see from \cref{fig:halfplane_a} that \cref{alg:eTFIRKA} places its poles along the branch cut of $\calA_G$ forming clusters towards its branch points at $-1$ and $1$. 
This is a known and recurrent phenomenon in rational approximation (see, e.g., \cite[Section 6.2]{NaSeTr18}, \cite{Sta97} or \cite{TrNaWe21}).

To conclude the experiment, we show the $\etwod$ relative error norm defined as 
\[
\frac{\|G-F_r\|_{\etwod}}{\|G\|_{\etwod}},
\]
for \cref{alg:eTFIRKA} and TF-IRKA at different values of $r$. 
Here, both algorithms provide approximately the same level of accuracy.

\subsection{BDF2 method}\label{sec:numbdf2}
For this example we apply the BDF2 method to the clamped beam LTI model from the SLICOT library (see \cite{slicotfom,ChaDoo05} and \cite[Section 6.9]{NaSeTr18}).
The considered transfer function is $G(z) = \bfc^\top\left(z\bfI -\bfA\right)^{-1} \bfb$ with $\bfA\in\R^{348\times 348}$ and $\bfc,\bfb\in\R^{348}$, and the conformal map $\varphi$ we use is defined in \eqref{eq:bdf2confmap} with $h=0.001$. 

In this experiment we compare the performance between \cref{alg:eTFIRKA} in the $\durenhardy$ setting and the AAA algorithm introduced in \cite{NaSeTr18}. 
The AAA approximants have the following rational barycentric formulation
\[
R(z) = \dfrac{\sum_{j=1}^{r}\frac{w_jg_j}{z-z_j}}{\sum_{j=1}^{r}\frac{w_j}{z-z_j}},
\]
with $\{z_{j}\}_{j=1}^r\in\C$ the support points, $\{w_{j}\}_{j=1}^r\in\C$ the weights, and  $\{g_{j}\}_{j=1}^r\in\C$ the data values. 
In the comparison of the two algorithms, we use AAA to construct three approximants:
\begin{itemize}
    \item $R_\unitdisk$: this results from approximating $G\circ\varphi$ by using points on the unit circle.
    This method is comparable to \cref{alg:eTFIRKA}, where $G$ is first composed with $\varphi$ and then approximated with a rational function. 
    \item $R_\A$: this results from approximating directly $G$ by using points in $\partial\A$.
    \item $R_{\unitdisk}^{\texttt{L}}$: this is the result of applying AAA to $G\circ\varphi$ with no limit on $r$ and with the \textit{Lawson phase} on (see \cite{NaTr20}).
\end{itemize}
While for $R_{\unitdisk}^{\texttt{L}}$ there is no restriction on $r$, for $R_{\unitdisk}$ and $R_{\A}$ we restrict the order of the approximant and turn off the Lawson phase. 
For this numerical test, we choose the sample points for AAA as follows. Let $\{\theta_j\}_{j=1}^{N}$, with $N=500$, be logarithmically spaced points in $[10^{-7},1]$, then 
\[
Z_1 = \begin{bmatrix} e^{\pi\iunit \theta_1},\dots,e^{\pi\iunit \theta_N},e^{-\pi\iunit \theta_1},\dots,e^{-\pi\iunit \theta_N} \end{bmatrix} \in \C^{1\times1000},
\textnormal{ and } 
Z_2 = \varphi(Z_1) \in \C^{1\times1000},
\] 
where $Z_1$ is used for the computation of $R_\unitdisk$ and $R_\unitdisk^{\texttt{L}}$, while $Z_2$ is for $R_{\A}$.
We chose $Z_1$ and $Z_2$ in this way to provide AAA with enough data points to characterize the effect of the poles of $G$ near the boundary $\partial\A$. 
As a matter of fact, the spectrum of $G$ is close to a section of $\partial\A$ that interesects the origin of the complex plane.
This translates to $Z_1$ having logarithmically spaced points that accumulate in 1. 
 
\begin{figure}[hbt]
    \centering
    \input{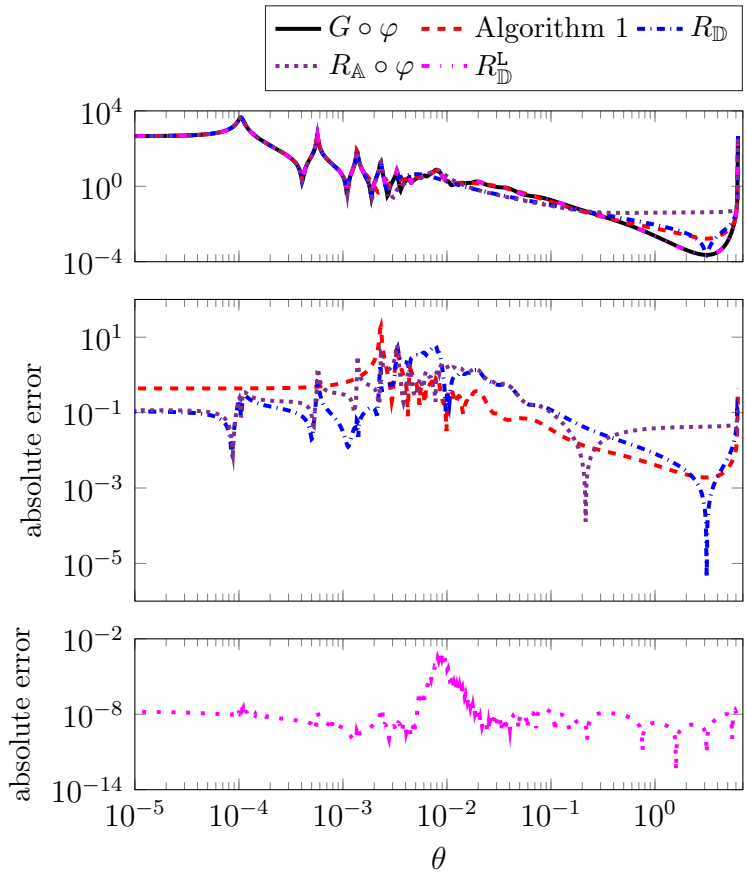}
    \caption{Evaluation of different approximants for $r=10$ and $h=0.001$.}
    \label{fig:BDF2evaldisk}
\end{figure}

In \cref{fig:BDF2evaldisk} we show the values and absolute error on the unit circle of the approximants computed with \cref{alg:eTFIRKA} and AAA for $r=10$. 
In more detail, the absolute error is defined as
\[
|G\circ\varphi(e^{\iunit\theta})-F_r\circ\varphi(e^{\iunit\theta})|,
\]
where $F_r$ is either $G_r$, as in \eqref{eq:Gr2}, $R_\unitdisk\circ\varphi^{-1}$, $R_\A$, or $R_\unitdisk^{\texttt{L}}\circ\varphi^{-1}$.  
Due to $R_\unitdisk^\texttt{L}$ having a greater advantage with respect to the other approximants, we divided the absolute errors in two plots. 
It is indeed possible to see in \cref{fig:BDF2evaldisk} that $R_\unitdisk^\texttt{L}$ reaches the lowest error. 
However, this result is provided with $r=99$ while the other approximants are constrained to $r=10$.
Let us then look at the middle plot of \cref{fig:BDF2evaldisk} where the absolute errors of $G_r\circ\varphi$, $R_\unitdisk$, and $R_\A\circ\varphi$ are shown. 
In the values of $\theta$ close to 0 until $\theta\approx 2\cdot10^{-3}$ both AAA approximations reach better results than \cref{alg:eTFIRKA}. 
This can be due to the accumulation of sample points of $Z_1$ and $Z_2$ around $\theta=0$. 
For the rest of the unit circle, instead, \cref{alg:eTFIRKA} provides an overall lower error. 
In \cref{fig:H2ABDF2} we show the resulting relative $\durenhardy$ error norm 
$$\frac{\|G-F_r\|_{\durenhardy}}{\|G\|_{\durenhardy}},$$
for different values of $r$.  
In \cref{fig:H2ABDF2} we omitted $R_\unitdisk^{\texttt{L}}$ due to its constant reduced order of $r=99$.
Nevertheless, it reaches a $\durenhardy$ relative error of $9.67\times 10^{-7}$, outperforming all the other approximants. 
On the other side, it is possible to see that \cref{alg:eTFIRKA} achieves lower $\durenhardy$ errors than $R_\unitdisk$ and $R_\A$ on most of the values of $r$.
\begin{figure}[hbt]
    \centering
    \begin{subfigure}[b]{0.49\textwidth}
        \centering
        \input{figures/H2AerrorBDF2.tex}
        \caption{}
        \label{fig:H2ABDF2}
    \end{subfigure}
    \begin{subfigure}[b]{0.49\textwidth}
        \centering
        \input{figures/errorBDF2.tex}
        \caption{}
        \label{fig:trajBDF2}
    \end{subfigure}    
    \caption{(Left)  Relative $\durenhardy$ error norm for \cref{alg:eTFIRKA} and AAA at different values of $r$ with $h=0.001$.
    (Right) Output trajectories and error between $y_k$ and $\widehat{y}_k$.}
    \label{fig:trajH2ABDF2}
\end{figure}
The reason why \cref{alg:eTFIRKA} does not always compute the best approximation can be due to the choice of the initial shifts which leads to potentially worse locally optimal reduced models. 
While \cref{fig:H2ABDF2} provides an interesting insight on the performance of \cref{alg:eTFIRKA}, we remind that the latter is designed to minimize the $\durenhardy$ norm and AAA is not.

We mentioned in \cref{sec:mor_disc} that applying \cref{alg:eTFIRKA} can be recasted as applying TF-IRKA to the discretized in time version of $G$. 
In addition, we provided an $\ell_\infty$ error bound between the output of the discretized system equivalent to $G\circ\varphi(z^{-1})$ and the ROM resulting from \cref{alg:eTFIRKA}. 
On this note, in \cref{fig:trajBDF2} we show the output error $|y_k-\widehat{y}_k|$ between \eqref{eq:ssbdf2} and the reduced model resulting from $G_r\circ\varphi(z^{-1})$, see e.g. \eqref{eq:ssGr_disc}, where $G_r$ is computed  $r=10$.
The responses $y_k$ and $\widehat{y}_k$ are a result of applying the discrete-time unit impulse of the form 
\begin{equation}\label{eq:discimpulse}
    u_k = \delta_k(m) := \begin{cases}
    &1, \; k=m,\\
    &0, \; k\neq m,
\end{cases}
\end{equation}
for $m=10$. 
It is particularly interesting to see that the absolute error respects the $\durenhardy$ boundary derived in \eqref{eq:linftybound}.
To note that for $u_k = \delta_k(m)$ we get $U(z) = z^{-m}$ and so $\|U\circ e^{\iunit\cdot}\|_{\mathcal{L}_2(0,2\pi)}=1$.

\subsection{BDF4 method}\label{sec:bdf4num}
\begin{figure}[hbt]
    \centering
        \includegraphics[width=0.9\textwidth]{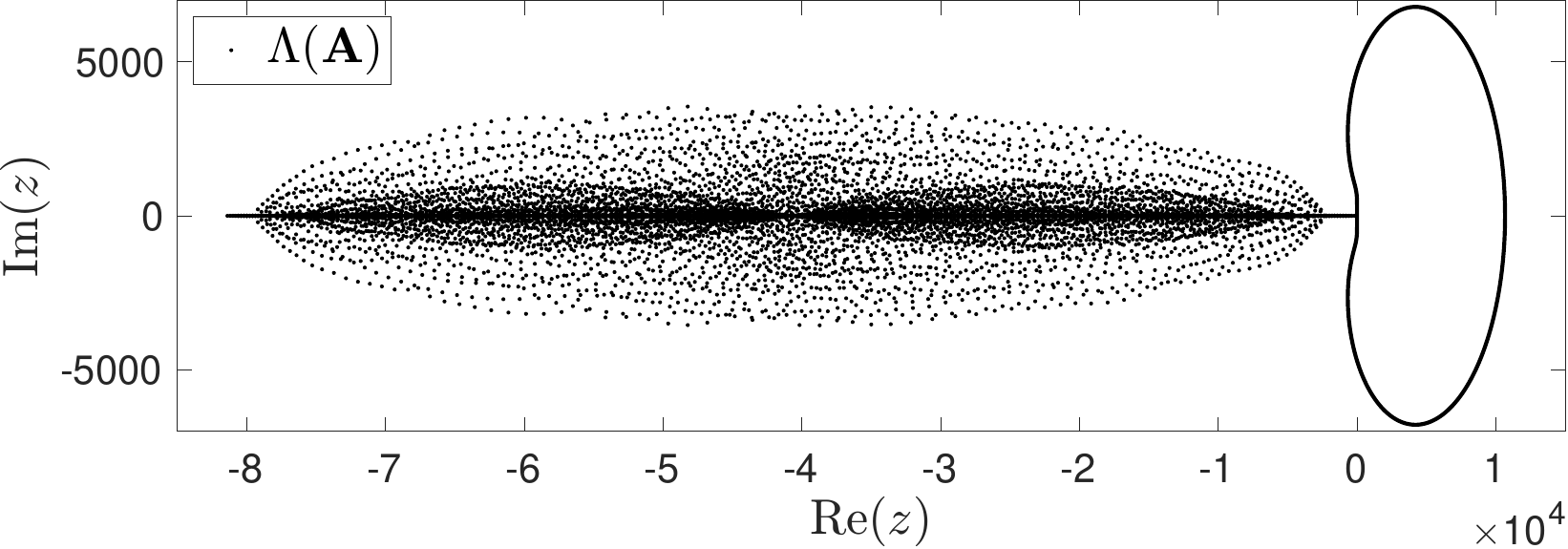}
    \caption{The spectrum of $G$ and the image $\A$ of the conformal map resulting from BDF4.} 
    \label{fig:eigsBDFinA}
\end{figure}

In this section we look at the BDF4 discretization technique discussed in \cref{sec:bdf}. 
For this experiment we consider a modified version of the finite-difference space discretization of a two-dimensional convection-diffusion equation on a unit square provided in \cite[Section C.1.1]{Pen99}.
More in detail, let $\Omega=(0,1)\times(0,1)$, we consider 
\begin{equation}\label{eq:convdiff}
    \begin{cases}
        w_t = \Delta w -  100(xw_x + y w_y) + \mathcal{I}_{[0.1,0.3]} u, &\text{ for } (t,x,y)\in(0,T]\times\Omega,\\
        w = 0, &\text{ for } (t,x,y)\in(0,T]\times\partial\Omega \\
        y = \mathcal{I}_{[0.7,0.9]}w, &\text{ for } t\in(0,T],
    \end{cases}
\end{equation} 
where $\mathcal{I}_{[a,b]}$ is the indicator function for the interval $[a,b]$ on the $x$ dimension in space.
After a spatial discretization with 100 grid points for each space dimension, we reformulate \eqref{eq:convdiff} as an LTI system with transfer function $G(z) = \bfc^\top \left(z\bfI - \bfA\right)^{-1}\bfb$, where $\bfA\in\R^{10000\times 10000}$, and $\bfc,\bfb\in\R^{10000}$.
The spectrum of $G$ is showed in \cref{fig:eigsBDFinA} along with the image $\A$ of $\varphi$ in \eqref{eq:confmapbdf4} for $h=0.001$.
It can be seen that $G$ is analytic in $\A$ making it an element of the $\durenhardy$ space.

\begin{figure}[hbt]
    \centering
    \begin{subfigure}[b]{0.5\textwidth}
        \input{figures/errorBDF.tex}
        \caption{}
        \label{fig:trajBDF}
    \end{subfigure}
    \begin{subfigure}[b]{0.45\textwidth}
        \input{figures/H2errorBDF.tex}
        \caption{}
        \label{fig:H2ABDF}
    \end{subfigure}
        
    \caption{(Left) Output trajectories and absolute error between $y_k$ and $\widehat{y}_k$. (Right) Relative $\durenhardy$ error norm for \cref{alg:eTFIRKA} at different values of $r$ with $h=0.001$.}
    \label{fig:trajH2ABDF}
\end{figure}
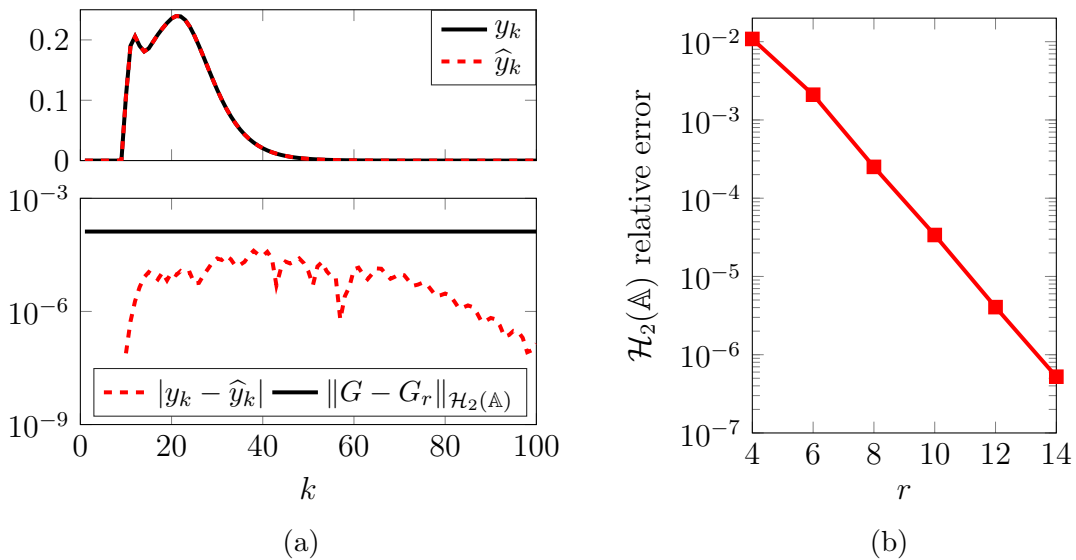

Similar to \cref{sec:numbdf2}, in \cref{fig:trajBDF} we show the impulse response of \eqref{eq:ssdiscbdf4} and \eqref{eq:ssGr_disc} resulting from the input \eqref{eq:discimpulse} for $r=10$.
Once again we can see that the error trajectory $|y_k - \widehat{y}_k|$ is bounded by the $\durenhardy$ error norm as proved in \eqref{eq:linftybound}.
Lastly, \cref{fig:H2ABDF} shows the convergence of \cref{alg:eTFIRKA} with respect to the $\durenhardy$ norm.

\paragraph*{Funding} The work of the authors was funded by the Deutsche
Forschungsgemeinschaft (DFG, German Research Foundation) - 504768428

\paragraph*{Conflict of interest} The authors declare no competing interests.

\begin{appendices}
\section{Proof \cref{th:E2Dintcond}} \label{proof:e2dintcond}
\begin{proof} 
    We consider a perturbed version $G_r^\varepsilon$ of the local minimizer $G_r$. 
    We can then write 
    \begin{align*}
        \|G-G_r\|_{\etwod}^2&\leq \|G-G_r^\varepsilon\|_{\etwod}^2 = \|G-G_r + G_r -G_r^\varepsilon\|_{\etwod}^2\\
        &= \|G-G_r\|^2_{\etwod} + 2\re\left\{\langle G-G_r,G_r-G_r^\varepsilon \rangle_{\etwod}\right\} + \|G_r-G_r^\varepsilon\|_{\etwod}^2,
    \end{align*}
    which then results in 
    \begin{equation}\label{eq:finalineq}
        0\leq 2\re\left\{\langle G-G_r,G_r-G_r^\varepsilon \rangle_{\etwod}\right\} + \|G_r-G_r^\varepsilon\|_{\etwod}^2.
    \end{equation}
    We then divide the proof in three sections where we perturb $d_r$, the residues $\phi_j$, and the poles $\lambda_j$.
    In all cases we consider an additive perturbation $\varepsilon e^{\iunit\theta}$ where $\varepsilon\in\R_{\geq0}$ and $\theta\in[0,2\pi]$.
    Additionally, for simplifying the nomenclature, we set $E:=G-G_r$.
    \paragraph{Perturbation of $d_r$:} Let
    \[
    G_r^\varepsilon (z) = \left(\sum_{j=1}^r \frac{\phi_j}{(\varphi^{-1}(z) - \lambda_j)} + (d_r +\varepsilon e^{\iunit\theta})\right)\frac{1}{[\varphi'(\varphi^{-1}(z))]^{1/2}}.
    \]
    We then have
    \begin{equation}\label{eq:d1}
        \|G_r-G_r^\varepsilon\|_{\etwod}^2 = \left\|-\varepsilon e^{\iunit\theta}\right\|_{\hardy}^2=\mathcal{O}(\varepsilon^2) \quad \text{for }\varepsilon\rightarrow 0.
    \end{equation}
    and
    \begin{equation}\label{eq:d2}
        \begin{aligned}
            \langle E,G_r-G_r^\varepsilon \rangle_{\etwod} &= \left\langle \calA_E,-\varepsilon e^{\iunit\theta} \right\rangle_{\hardy} = -\varepsilon e^{-\iunit\theta}\calA_E\left(0\right).
        \end{aligned}
    \end{equation}
    Replacing the two terms in \eqref{eq:finalineq} with \eqref{eq:d1} and \eqref{eq:d2} we get 
    \begin{equation*}
        2\re\left\{e^{-\iunit\theta}\calA_E\left(0\right)\right\}\leq   \mathcal{O}(\varepsilon).
    \end{equation*}
    Choosing $\theta = \textsf{arg}\left(\calA_E\left(0\right)\right)$ leads to 
    \begin{equation*}
        0\leq2\left|\calA_E\left(0\right)\right|\leq   \mathcal{O}(\varepsilon),
    \end{equation*}
    which for $\varepsilon\rightarrow 0$ results in 
    \begin{equation*}
        G\circ\varphi\left(0\right) = G_r\circ\varphi\left(0\right),
    \end{equation*}
    proving the first equality in \eqref{eq:optintE2D}.

    \paragraph{Perturbation of the residue:} Consider
    \[
    G_r^\varepsilon (z) = \left(\sum_{j\neq p}^r \frac{\phi_j}{(\varphi^{-1}(z) - \lambda_j)} +\frac{\phi_p+\varepsilon e^{\iunit\theta}}{(\varphi^{-1}(z) - \lambda_p)} + d_r\right)\frac{1}{[\varphi'(\varphi^{-1}(z))]^{1/2}}.
    \]
    We then have 
    \begin{equation}\label{eq:res1}
        \|G_r-G_r^\varepsilon\|_{\etwod}^2 = \left\|\frac{-\varepsilon e^{\iunit\theta}}{\cdot - \lambda_p}\right\|_{\hardy}^2=\mathcal{O}(\varepsilon^2) \quad \text{for }\varepsilon\rightarrow 0,
    \end{equation}
    and 
    \begin{equation}\label{eq:res2}
        \begin{aligned}
            \langle E,G_r-G_r^\varepsilon \rangle_{\etwod} &= \left\langle \calA_E,\frac{-\varepsilon e^{\iunit\theta}}{\cdot - \lambda_p} \right\rangle_{\hardy} = \frac{\varepsilon e^{-\iunit\theta}}{\overline{\lambda_p}}\calA_E\left(\frac{1}{\overline{\lambda_p}}\right),
        \end{aligned}
    \end{equation}
    where we used \eqref{eq:evalinnerproductE2D}.
    Replacing the terms in \eqref{eq:finalineq} with \eqref{eq:res1} and \eqref{eq:res2} we get 
    \begin{equation*}
        2\re\left\{\frac{-e^{-\iunit\theta}}{\overline{\lambda_p}}\calA_E\left(\frac{1}{\overline{\lambda_p}}\right)\right\}\leq   \mathcal{O}(\varepsilon).
    \end{equation*}
    Choosing $\theta = \textsf{arg}\left(-\frac{1}{\overline{\lambda_p}}\calA_E\left(\frac{1}{\overline{\lambda_p}}\right)\right)$ leads to 
    \begin{equation*}
        0\leq2\left|\frac{1}{\overline{\lambda_p}}\calA_E\left(\frac{1}{\overline{\lambda_p}}\right)\right|\leq   \mathcal{O}(\varepsilon),
    \end{equation*}
    which for $\varepsilon\rightarrow 0$ results in 
    \begin{equation}\label{eq:int1}
        G\circ\varphi\left(\frac{1}{\overline{\lambda}_p}\right) = G_r\circ\varphi\left(\frac{1}{\overline{\lambda}_p}\right),
    \end{equation}
    proving the second equality in \eqref{eq:optintE2D}.
    
    \paragraph{Perturbation of the pole:} Let
    \[
    G_r^\varepsilon (z) = \left(\sum_{j\neq p}^r \frac{\phi_j}{(\varphi^{-1}(z) - \lambda_j)} +\frac{\phi_p}{(\varphi^{-1}(z) - \lambda_p-\varepsilon e^{\iunit\theta})} + d_r\right)\frac{1}{[\varphi'(\varphi^{-1}(z))]^{1/2}}.
    \]
    we then have
    \begin{equation}\label{eq:pol1}
        \left\|G_r-G_r^\varepsilon\right\|_{\etwod}^2 = \left\|\frac{\phi_p}{\cdot - \lambda_p} - \frac{\phi_p}{\cdot - (\lambda_p + \varepsilon e^{\iunit\theta})}\right\|_{\hardy}^2 = \mathcal{O}(\varepsilon^2) \quad \text{for }\varepsilon\rightarrow 0,
    \end{equation}
    and
    \begin{align*}
        \langle E,G_r-G_r^\varepsilon \rangle_{\etwod} &= \left\langle \calA_E,\frac{\phi_p}{\cdot - \lambda_p} - \frac{\phi_p}{\cdot - (\lambda_p + \varepsilon e^{\iunit\theta})}\right\rangle_{\hardy}\\        
        & = -\frac{\overline{\phi_p}}{\overline{\lambda_p}}\calA_E\left(\frac{1}{\overline{\lambda_p}}\right) + \frac{\overline{\phi_p}}{\overline{\lambda_p}+\varepsilon e^{-\iunit\theta}}\calA_E\left(\frac{1}{\overline{\lambda_p} + \varepsilon e^{-\iunit\theta}}\right)
    \end{align*}
    Given \eqref{eq:int1} we have that $\calA_E\left(1/\overline{\lambda_p}\right)=0$ leading to 
    \[
        \langle E,G_r-G_r^\varepsilon \rangle_{\etwod} =\frac{\overline{\phi_p}}{\overline{\lambda_p}+\varepsilon e^{-\iunit\theta}}\calA_E\left(\frac{1}{\overline{\lambda_p} + \varepsilon e^{-\iunit\theta}}\right).
    \]
    To note that $\calA_E\in \hardy$ is analytic around $1/\overline{\lambda_p}$ as $|\lambda_p|>1$.
    We then consider the Taylor expansion of $\calA_E(1/z)$ around $\overline{\lambda_p}$ for $\varepsilon\rightarrow0$. This results in 
    \begin{equation}\label{eq:pol2}
        \begin{aligned}
            \langle E,G_r-G_r^\varepsilon \rangle_{\etwod} &= \frac{\overline{\phi_p}}{\overline{\lambda_p} + \varepsilon e^{-\iunit\theta}}\left(\calA_E\left(\frac{1}{\overline{\lambda_p}}\right) -\varepsilon e^{-\iunit\theta}\calA_E'\left(\frac{1}{\overline{\lambda_p}}\right)\frac{1}{\overline{\lambda_p}^2} + \mathcal{O}(\varepsilon^2)\right)\\
            &= \frac{\overline{\phi_p}}{\overline{\lambda_p}+\varepsilon e^{-\iunit\theta}}\left(-\varepsilon e^{-\iunit\theta}\calA_E'\left(\frac{1}{\overline{\lambda_p}}\right)\frac{1}{\overline{\lambda_p}^2} + \mathcal{O}(\varepsilon^2)\right).
        \end{aligned}
    \end{equation}
    Again, substituting the terms in \eqref{eq:finalineq} with \eqref{eq:pol1} and \eqref{eq:pol2} results in 
    \begin{align*}
        2\re\left\{\frac{\overline{\phi_p}e^{-\iunit\theta}}{\overline{\lambda_p}+\varepsilon e^{-\iunit\theta}}\calA_E'\left(\frac{1}{\overline{\lambda_p}}\right)\frac{1}{\overline{\lambda_p}^2}\right\}\leq\mathcal{O}(\varepsilon)
    \end{align*}
    For $\varepsilon\rightarrow0$ and choosing $\theta=\textsf{arg}\left(\frac{\overline{\phi_p}}{\overline{\lambda_p}}\calA_E'\left(\frac{1}{\overline{\lambda_p}}\right)\frac{1}{\overline{\lambda_p}^2}\right)$ we then get 
    \[
        \calA_E'\left(\frac{1}{\overline{\lambda_p}}\right) = 0,
    \]
    which can be expanded into
    \[
    E'\circ\varphi\left(\frac{1}{\overline{\lambda_p}}\right)\varphi'\left(\frac{1}{\overline{\lambda_p}}\right)\left[\varphi'\left(\frac{1}{\overline{\lambda_p}}\right)\right]^{1/2} + \underbrace{E\circ\varphi\left(\frac{1}{\overline{\lambda_p}}\right)}_{=0 \text{ due to \eqref{eq:int1}}}\frac{\mathrm{d}}{\mathrm{d}z}\left[\varphi'\left(z\right)\right]^{1/2}\big|_{z=\frac{1}{\overline{\lambda_p}}} = 0.
    \]
    Due to $\varphi'$ not vanishing in $1/\overline{\lambda_p}$ for $|\lambda_p|>1$ we get 
    \[
        G'\circ\varphi\left(\frac{1}{\overline{\lambda}_p}\right) = G_r'\circ\varphi\left(\frac{1}{\overline{\lambda}_p}\right),
    \]
    proving the last equality in \eqref{eq:optintE2D}.
\end{proof}

\end{appendices}

\bibliographystyle{siam} 
\bibliography{refs}
		
\end{document}

%% file: pkgssetup.tex
\usepackage[utf8]{inputenc}

\usepackage{amsfonts}
\usepackage{amsmath}
\usepackage{amssymb}
\usepackage{amsthm}
\usepackage{algpseudocode}
\usepackage{algorithm}%
\usepackage{graphicx} 
\usepackage{pgfgantt}
\usepackage[symbol]{footmisc}
\usepackage{pdflscape}
\usepackage{appendix}
\usepackage{pgfplots}
\usepackage{tikz}
    \pgfplotsset{compat=newest}
    \usetikzlibrary{plotmarks}
    \usetikzlibrary{calc}
    \usetikzlibrary{positioning}
    \usetikzlibrary{arrows.meta}
    \usetikzlibrary{matrix}
    \usepgfplotslibrary{patchplots}
	\usetikzlibrary{intersections,decorations.markings}	
\usepackage[framemethod=TikZ]{mdframed}
\usepackage{xspace}
\usepackage{hyperref}
\usepackage{cleveref}
\usepackage{url}
\usepackage{wrapfig}
\usepackage{subcaption}
\usepackage{marginnote}
\usepackage{titlepic}
\usepackage{caption}
\usepackage{subcaption}
\graphicspath{ {./figures/} }

\usetikzlibrary{shapes.misc}
\tikzset{cross/.style={cross out, draw=black, fill=none, minimum size=2*(#1-\pgflinewidth), inner sep=0pt, outer sep=0pt}, cross/.default={2pt}}

\newtheorem{lemma}{Lemma}
\newtheorem{definition}{Definition}

\newtheorem{theorem}{Theorem}
\newtheorem{corollary}{Corollary}
\newtheorem{assumption}{Assumption}
\newtheorem{remark}{Remark}

\newcommand\C{\ensuremath{\mathbb{C}}}
\newcommand\R{\ensuremath{\mathbb{R}}}
\newcommand\A{\ensuremath{\mathbb{A}}}
\newcommand\unitdisk{\ensuremath{\mathbb{D}}}

\newcommand\hardy{\ensuremath{\mathcal{H}_2}}
\newcommand\durenhardy{\ensuremath{\mathcal{H}_2(\mathbb{A})}}
\newcommand\etwod{\ensuremath{\mathcal{E}_2(\mathbb{A})}}
\newcommand\re{\ensuremath{\text{Re}}}
\newcommand\im{\ensuremath{\text{Im}}}

\newcommand\iunit{\ensuremath{\mathrm{i}}}
\newcommand\bfI{\ensuremath{\mathbf{I}}}
\newcommand\bfc{\ensuremath{\mathbf{c}}}
\newcommand\bfb{\ensuremath{\mathbf{b}}}

\newcommand\bfA{\ensuremath{\mathbf{A}}}
\newcommand\bfE{\ensuremath{\mathbf{E}}}

\newcommand\calA{\ensuremath{\mathcal{A}}}

%% file: figures/iR_E2Cplus.tex
%
%
\begin{tikzpicture}

\begin{axis}[%
width=5cm,
height=4cm,
at={(0,0)},
scale only axis,
xmode=log,
xmin=0.01,
xmax=1000,
xminorticks=true,
xlabel style={font=\color{white!15!black}},
xlabel={$\iunit\omega$},
ymode=log,
ymin=1e-06,
ymax=0.01,
yminorticks=true,
ylabel style={font=\color{white!15!black}},
ylabel={absolute error},
axis background/.style={fill=white},
legend style={at={(1,1)},legend cell align=left, align=left, draw=white!15!black}
]
\addplot [color=red, line width = 1.5pt]
  table[row sep=crcr]{%
0.01	0.00193726055734851\\
0.0146273335620113	0.00161573243933652\\
0.0195114834684662	0.00140009280267241\\
0.0327729484992338	0.00107394363374601\\
0.0484937406733523	0.000881928169077852\\
0.14831025143361	0.000505150732147568\\
0.381576466127125	0.000314994648086638\\
1.48652484499786	0.000158870971782887\\
2.43998629725955	0.000123095020499647\\
3.40825854742345	0.000103080223836712\\
4.49420266211914	8.85075166049333e-05\\
5.79112264764176	7.6522650975273e-05\\
7.46230289139111	6.57666183333842e-05\\
10.423606739764	5.35145443557212e-05\\
14.728827239075	4.32994610314456e-05\\
17.916503273639	3.86201055240881e-05\\
21.2974853574552	3.51093548837488e-05\\
25.0264009641792	3.23082312920936e-05\\
29.0712337727258	3.00784801386863e-05\\
33.3828586473176	2.83081960702117e-05\\
38.3339510176661	2.67924393156798e-05\\
43.5149650092505	2.56148240281435e-05\\
48.8302208687789	2.47182058979464e-05\\
54.1668691103316	2.40586790707787e-05\\
59.3982669392036	2.35995938130737e-05\\
65.1349094627281	2.32734910499978e-05\\
70.6071771413778	2.31095397181464e-05\\
76.5391938823016	2.30699373208188e-05\\
82.0188949920221	2.31411684988424e-05\\
88.9096598952917	2.33495918660161e-05\\
96.379347996158	2.36879021803948e-05\\
105.687597118481	2.42121987574499e-05\\
122.769104798836	2.5255675496866e-05\\
139.361927422414	2.6122032267836e-05\\
151.070330448666	2.65512708949021e-05\\
161.88596901782	2.67863852018385e-05\\
171.488196987054	2.68605710553541e-05\\
181.659978837533	2.68044921133585e-05\\
192.435097523033	2.66035985479171e-05\\
203.849339825246	2.62478697199492e-05\\
213.466303332425	2.58483950705385e-05\\
223.53696459098	2.53468160003434e-05\\
234.082727617829	2.4745738839302e-05\\
247.967289250216	2.38624111965884e-05\\
262.675410372384	2.2846792630789e-05\\
278.255940220713	2.17192990965312e-05\\
294.760625512486	2.05048512273027e-05\\
312.244282309286	1.92317378937917e-05\\
334.598912054997	1.76691723656175e-05\\
362.710025233065	1.5864250315878e-05\\
402.350554886929	1.36868479119244e-05\\
483.820966492596	1.05124379138032e-05\\
530.548052536957	9.29666488642324e-06\\
581.788007434494	8.28986904059294e-06\\
645.37154016467	7.34657355051081e-06\\
741.088151564157	6.30386227133317e-06\\
933.189771573324	4.92550141196414e-06\\
1000	4.57958597975093e-06\\
};
\addlegendentry{Algorithm 1}

\addplot [color=blue, dashed, line width = 1.5pt]
  table[row sep=crcr]{%
0.01	0.00487783106767853\\
0.0123052400435926	0.00435872697622021\\
0.0151418932530435	0.00391789085317472\\
0.0190669084051225	0.00350089049119314\\
0.0468458011587305	0.00227247721435915\\
0.0701206358900718	0.0018495583064842\\
0.10991097009295	0.00147255289168545\\
0.845136633068472	0.000530434873753725\\
1.35560178532937	0.000416287778307854\\
1.93770333747799	0.000344879467163884\\
2.49687842888433	0.000300195092623047\\
3.072468842709	0.00026651984434743\\
3.65226736430818	0.000240087405593986\\
4.29173237842216	0.000216551494682081\\
4.92824957004051	0.000197142462465924\\
5.65917016324624	0.000178382244409334\\
6.42403365939419	0.000161746887062344\\
7.20871503378214	0.00014711588753727\\
8.08924348680594	0.000133004657590353\\
9.07732652521022	0.000119483849843932\\
10.1861017015598	0.000106636695416209\\
11.4303112911448	9.45489495214009e-05\\
12.9751716865759	8.22220777429915e-05\\
14.728827239075	7.09878223980407e-05\\
16.9132951702965	6.00311561528413e-05\\
19.8745954958099	4.89678893432822e-05\\
24.1759407916913	3.79150831841002e-05\\
42.5234633452869	1.80299477874131e-05\\
49.3962174387833	1.49722129181269e-05\\
56.0723488285204	1.28942665733911e-05\\
62.9214610961035	1.13449612804316e-05\\
69.7981390783067	1.01874857975516e-05\\
76.5391938823016	9.32315609981358e-06\\
82.9695852083491	8.68006018377643e-06\\
89.9402217409205	8.13239608393698e-06\\
97.496491834841	7.67144878111081e-06\\
104.476597156081	7.33917033264568e-06\\
111.956431948388	7.05950881298102e-06\\
119.971773543589	6.82720497511643e-06\\
128.56096069433	6.63665229073628e-06\\
139.361927422414	6.45897355708166e-06\\
151.070330448666	6.31927767440422e-06\\
167.580786453077	6.17695005382416e-06\\
211.02034285686	5.88689342377986e-06\\
228.74908173557	5.75495810632552e-06\\
245.126006203334	5.61685740702649e-06\\
262.675410372384	5.45111671977307e-06\\
278.255940220713	5.29000132491658e-06\\
294.760625512486	5.10747532708514e-06\\
312.244282309286	4.90414068959225e-06\\
330.764978074424	4.68164697191572e-06\\
350.384224529068	4.44265832393905e-06\\
371.167181947577	4.19075743392055e-06\\
397.740302405804	3.8775804396752e-06\\
426.215882901532	3.56052763302141e-06\\
462.024137175131	3.19745414779274e-06\\
518.459354389291	2.71562895910408e-06\\
623.440188862786	2.08850123157496e-06\\
691.575882873852	1.81901796502152e-06\\
767.15811767793	1.59856111934881e-06\\
860.864769614924	1.39625279363405e-06\\
1000	1.18015388643607e-06\\
};
\addlegendentry{TF-IRKA}

\end{axis}
\end{tikzpicture}%

%% file: figures/partialD_E2Cplus.tex
%
%
\begin{tikzpicture}

\begin{axis}[%
  width=5cm,
  height=4cm,
  at={(0,0)},
scale only axis,
xmin=0,
xmax=6.28,
xlabel style={font=\color{white!15!black}},
xlabel={$\theta$},
xtick={0, 3.14, 6.28}, 
xticklabels={$0$, $\pi$, $2\pi$},
ymode=log,
ymin=0.00015,
ymax=0.1,
yminorticks=true,
ylabel style={font=\color{white!15!black}},
ylabel={absolute error},
axis background/.style={fill=white},
legend style={at={(1,1)},legend cell align=left, align=left, draw=white!15!black}
]
\addplot [color=red, line width = 1.5pt]
  table[row sep=crcr]{%
0.00628947478196196	0.00423208592183918\\
0.012578949563923	0.00300623068736229\\
0.018868424345885	0.00181614495472596\\
0.0251578991278461	0.00129837990830036\\
0.031447373909808	0.00104986816635032\\
0.0377368486917691	0.000906437103275277\\
0.044026323473731	0.000811529125880622\\
0.0503157982556921	0.000743146580565888\\
0.0566052730376541	0.000691125187375201\\
0.0628947478196151	0.000650036690986717\\
0.0691842226015771	0.000616664663701822\\
0.0754736973835382	0.000588959526414732\\
0.0817631721655001	0.000565544575667097\\
0.0880526469474621	0.000545457196544887\\
0.100631596511385	0.000512665460556603\\
0.113210546075308	0.000486873521869475\\
0.125789495639231	0.000465893468001803\\
0.138368445203154	0.00044835253833153\\
0.150947394767077	0.000433348599418437\\
0.163526344331	0.000420267287164695\\
0.182394768676885	0.000403343109693212\\
0.201263193022769	0.000388824297285487\\
0.220131617368654	0.000376090469914972\\
0.2452895164965	0.000361202959337838\\
0.270447415624346	0.00034814703093781\\
0.295605314752192	0.000336522073461809\\
0.327052688662	0.000323600560351695\\
0.358500062571808	0.000312136730731181\\
0.389947436481616	0.000301877326143841\\
0.427684285173386	0.00029089510590188\\
0.465421133865155	0.000281133840785604\\
0.503157982556924	0.000272405416540458\\
0.547184306030655	0.000263331508061853\\
0.591210629504386	0.000255276241370486\\
0.635236952978117	0.00024809233922934\\
0.685552751233809	0.000240798857248616\\
0.735868549489501	0.000234349128165248\\
0.792473822527155	0.00022796265386437\\
0.849079095564809	0.000222374701517781\\
0.905684368602463	0.000217483648671583\\
0.968579116422078	0.000212767559421746\\
1.03147386424169	0.000208721029766952\\
1.09436861206131	0.000205271946915591\\
1.16355283466289	0.000202099307378491\\
1.23273705726446	0.000199518724276049\\
1.30192127986604	0.000197482797190623\\
1.37110550246762	0.000195955568669765\\
1.44657919985116	0.0001948392630655\\
1.52205289723469	0.000194275550741631\\
1.59752659461823	0.000194254040289619\\
1.67300029200177	0.000194775481626245\\
1.74847398938531	0.000195851909217955\\
1.81765821198689	0.000197346453275493\\
1.88684243458846	0.000199355616971332\\
1.95602665719004	0.00020191832872279\\
2.02521087979162	0.00020508608307545\\
2.08810562761123	0.000208546698677878\\
2.15100037543085	0.000212625858725167\\
2.21389512325046	0.000217403381411054\\
2.27050039628812	0.000222382946501432\\
2.32710566932577	0.000228101634659041\\
2.37742146758147	0.000233898606609038\\
2.42773726583716	0.000240473536702356\\
2.47805306409285	0.00024795627419039\\
2.52207938756658	0.0002553788097884\\
2.56610571104031	0.000263768750494347\\
2.60384255973208	0.000271871343682776\\
2.64157940842385	0.000280975782168041\\
2.67931625711562	0.000291281497684817\\
2.71076363102543	0.000300974903968038\\
2.74221100493523	0.000311883563732276\\
2.77365837884504	0.000324267102330981\\
2.79881627797289	0.000335462791391467\\
2.82397417710074	0.000348047751928629\\
2.84913207622858	0.000362314539422169\\
2.87428997535643	0.00037864433432107\\
2.89315839970231	0.000392543971280753\\
2.9120268240482	0.000408162823528354\\
2.93089524839408	0.000425857421119791\\
2.94976367273997	0.000446097146229872\\
2.96863209708585	0.000469520444811413\\
2.98750052143174	0.000497033404985124\\
3.00007947099566	0.000518290073056055\\
3.01265842055958	0.000542531022023545\\
3.0252373701235	0.000570572367213035\\
3.03781631968743	0.000603605755719957\\
3.05039526925135	0.000643438694535151\\
3.06297421881527	0.000692925024120353\\
3.06926369359724	0.00072265077891141\\
3.0755531683792	0.000756767947189508\\
3.08184264316116	0.000796410289465566\\
3.08813211794312	0.000843110627446474\\
3.09442159272508	0.000898988343462526\\
3.10071106750704	0.000967072320451899\\
3.10700054228901	0.00105192956543608\\
3.11329001707097	0.00116111812923464\\
3.11957949185293	0.00130937929334943\\
3.12586896663489	0.00153450056408566\\
3.13215844141685	0.00198198859540226\\
3.13844791619881	0.00371282414229187\\
3.14473739098077	0.00371282725923264\\
3.15102686576274	0.00198199147369728\\
3.1573163405447	0.00153450290390386\\
3.16360581532666	0.00130938127704222\\
3.16989529010862	0.00116111988480391\\
3.17618476489058	0.00105193116453829\\
3.18247423967254	0.0009670738035899\\
3.18876371445451	0.000898989734784561\\
3.19505318923647	0.000843111942295747\\
3.20134266401843	0.000796411538424174\\
3.20763213880039	0.000756769138032073\\
3.22021108836431	0.0006929261157822\\
3.23279003792824	0.000643439704087614\\
3.24536898749216	0.000603606696084133\\
3.25794793705608	0.000570573248656278\\
3.27052688662001	0.000542531852803356\\
3.28310583618393	0.000518290859899352\\
3.29568478574785	0.000497034153411383\\
3.31455321009374	0.000469521143823704\\
3.33342163443962	0.000446097803555999\\
3.35229005878551	0.000425858042701809\\
3.37115848313139	0.000408163413971671\\
3.39002690747727	0.000392544534266336\\
3.40889533182316	0.000378644872821984\\
3.43405323095101	0.000362315049017139\\
3.45921113007885	0.000348048236072506\\
3.4843690292067	0.000335463252891385\\
3.50952692833454	0.000324267543500846\\
3.54097430224435	0.00031188398223169\\
3.57242167615416	0.00030097530231071\\
3.60386905006397	0.000291281877999272\\
3.64160589875574	0.000280976143201775\\
3.67934274744751	0.000271871687590302\\
3.71707959613927	0.000263769079115237\\
3.76110591961301	0.000255379122508414\\
3.80513224308674	0.00024795657274624\\
3.84915856656047	0.000241356498217516\\
3.89947436481616	0.00023467583428945\\
3.94979016307185	0.000228787699242471\\
4.00639543610951	0.000222980214096005\\
4.06300070914716	0.000217923447194218\\
4.11960598218481	0.000213522837822013\\
4.18250073000443	0.000209310971764579\\
4.24539547782405	0.000205732190964168\\
4.30829022564366	0.000202719957341212\\
4.37747444824524	0.000199997917820216\\
4.44665867084681	0.000197842129009262\\
4.51584289344839	0.000196210450448977\\
4.58502711604997	0.000195071872896153\\
4.66050081343351	0.000194367338954157\\
4.73597451081705	0.000194208492296836\\
4.81144820820058	0.000194590863826978\\
4.88692190558412	0.000195521114108658\\
4.9561061281857	0.000196870305336387\\
5.02529035078728	0.000198717633778517\\
5.09447457338885	0.000201095313797704\\
5.16365879599043	0.000204046365186351\\
5.23284301859201	0.000207626586112846\\
5.29573776641162	0.000211487017332017\\
5.35863251423124	0.00021599367425593\\
5.42152726205085	0.000221230849081389\\
5.47813253508851	0.000226655929460603\\
5.53473780812616	0.000232854949030613\\
5.58505360638185	0.000239112246725706\\
5.63536940463755	0.00024618282331133\\
5.68568520289324	0.000254199745347319\\
5.72971152636697	0.000262122281713606\\
5.7737378498407	0.000271040291013758\\
5.81147469853247	0.00027961187753026\\
5.84921154722424	0.000289189176173452\\
5.88694839591601	0.000299952680386679\\
5.91839576982582	0.000309995297566094\\
5.94984314373562	0.000321200532180766\\
5.98129051764543	0.000333805859926744\\
6.00644841677328	0.000345118896981122\\
6.03160631590112	0.000357785860270746\\
6.05676421502897	0.000372167654516443\\
6.07563263937485	0.000384404868740308\\
6.09450106372074	0.000398272450270139\\
6.11336948806662	0.000414307762571682\\
6.12594843763055	0.000426597898694729\\
6.13852738719447	0.000440577042679297\\
6.15110633675839	0.000456757796860341\\
6.16368528632232	0.00047588170795544\\
6.17626423588624	0.000499056038039343\\
6.18884318545016	0.000528002390527359\\
6.19513266023212	0.000545457109091864\\
6.20142213501409	0.000565544486880289\\
6.20771160979605	0.000588959436521948\\
6.21400108457801	0.00061666457338381\\
6.22029055935997	0.000650036600297591\\
6.22658003414193	0.000691125096644437\\
6.23286950892389	0.000743146489927872\\
6.23915898370586	0.00081152903541842\\
6.24544845848782	0.00090643701338902\\
6.25173793326978	0.00104986807871105\\
6.25802740805174	0.0012983798321014\\
6.2643168828337	0.00181614491654362\\
6.27060635761566	0.003006230733587\\
6.27689583239763	0.00423208600352821\\
6.28318530717959	0.00316131721514194\\
};
\addlegendentry{Algorithm 1}

\addplot [color=blue, dashed, line width = 1.5pt]
  table[row sep=crcr]{%
0.00628947478196196	0.00108732678418712\\
0.012578949563923	0.00070213818253901\\
0.018868424345885	0.000545404848922414\\
0.0251578991278461	0.000506426811586041\\
0.031447373909808	0.000504446369191465\\
0.0377368486917691	0.000516001555048848\\
0.0503157982556921	0.000552625410274649\\
0.0628947478196151	0.000592522989083353\\
0.0754736973835382	0.000629564858876439\\
0.0880526469474621	0.000661911559116559\\
0.100631596511385	0.000689207497403488\\
0.113210546075308	0.000711665423511885\\
0.125789495639231	0.000729714445208909\\
0.138368445203154	0.000743852198609041\\
0.150947394767077	0.000754578672605711\\
0.163526344331	0.000762365724433597\\
0.176105293894923	0.000767643635516523\\
0.188684243458846	0.000770796410537973\\
0.201263193022769	0.00077216174090898\\
0.220131617368654	0.000771489889083263\\
0.239000041714539	0.000768274285001069\\
0.257868466060423	0.000763169324250013\\
0.283026365188269	0.000754293860370488\\
0.314473739098077	0.000741104820120904\\
0.358500062571808	0.000720780398407818\\
0.534605356466732	0.000642818341639621\\
0.597500104286347	0.000619639804162947\\
0.654105377324001	0.000601033427114501\\
0.717000125143617	0.000582674285079389\\
0.779894872963232	0.000566524443128838\\
0.842789620782847	0.00055235909961299\\
0.905684368602463	0.000539977521752314\\
0.968579116422078	0.000529206102611316\\
1.03776333902366	0.000519041630348273\\
1.10694756162523	0.000510482686982708\\
1.17613178422681	0.000503394243015164\\
1.24531600682839	0.000497666739239792\\
1.31450022942996	0.000493212905349433\\
1.3899739268135	0.000489728116701155\\
1.46544762419704	0.0004876164695697\\
1.54092132158058	0.000486839527597787\\
1.61639501896412	0.000487382790013357\\
1.69186871634766	0.000489256363348342\\
1.7673424137312	0.000492496830542039\\
1.83652663633277	0.000496724099489639\\
1.90571085893435	0.000502231861416298\\
1.97489508153593	0.000509125854265874\\
2.0440793041375	0.000517547413484108\\
2.10697405195712	0.00052668423673199\\
2.16986879977673	0.000537415954216957\\
2.22647407281439	0.000548625651769226\\
2.28307934585204	0.00056152064425855\\
2.3396846188897	0.000576362342589516\\
2.39000041714539	0.000591453163246648\\
2.44031621540108	0.000608629392919672\\
2.48434253887481	0.00062565054484228\\
2.52836886234854	0.000644843460403583\\
2.57239518582227	0.000666582339738314\\
2.61013203451404	0.000687598029993946\\
2.64786888320581	0.000711208776354114\\
2.68560573189758	0.000737896943914748\\
2.71705310580739	0.000762937176039231\\
2.7485004797172	0.000791026538109946\\
2.779947853627	0.000822793569290416\\
2.80510575275485	0.000851417033755946\\
2.8302636518827	0.000883517047898986\\
2.85542155101054	0.000919874016226897\\
2.88057945013839	0.000961554522974677\\
2.89944787448427	0.000997185147419852\\
2.91831629883016	0.0010374929488885\\
2.93718472317604	0.00108361590652848\\
2.95605314752193	0.00113709080787734\\
2.96863209708585	0.00117783634153822\\
2.98121104664977	0.00122358838418554\\
2.9937899962137	0.00127537956003291\\
3.00636894577762	0.00133452866775292\\
3.01894789534154	0.00140274559601458\\
3.03152684490547	0.00148229349818379\\
3.04410579446939	0.00157625884458511\\
3.05668474403331	0.00168904738485361\\
3.06926369359724	0.00182741338360605\\
3.08184264316116	0.00200292798048603\\
3.08813211794312	0.0021110464410972\\
3.09442159272508	0.00223897162156219\\
3.10071106750704	0.00239519973582339\\
3.10700054228901	0.00259481570477323\\
3.11329001707097	0.00286678747993676\\
3.11957949185293	0.00327207027539948\\
3.12586896663489	0.00395495850733279\\
3.13215844141685	0.00532730859302721\\
3.13844791619881	0.00923665765654602\\
3.14473739098077	0.00923665728385182\\
3.15102686576274	0.00532730812392214\\
3.1573163405447	0.003954958042506\\
3.16360581532666	0.00327206984128302\\
3.16989529010862	0.00286678708120305\\
3.17618476489058	0.00259481533863177\\
3.18247423967254	0.00239519939763137\\
3.18876371445451	0.00223897130665784\\
3.19505318923647	0.00211104614558019\\
3.20763213880039	0.00190948385149404\\
3.22021108836431	0.00175446924351047\\
3.23279003792824	0.00162995108572937\\
3.24536898749216	0.00152723759518134\\
3.25794793705608	0.00144093363170034\\
3.27052688662001	0.00136737743676556\\
3.28310583618393	0.0013039380230152\\
3.29568478574785	0.00124865421024206\\
3.30826373531177	0.00120002758212792\\
3.32713215965766	0.00113709064552107\\
3.34600058400354	0.00108361575223922\\
3.36486900834943	0.00103749280176453\\
3.38373743269531	0.000997185006819384\\
3.4026058570412	0.000961554388245285\\
3.42776375616904	0.00091987388849751\\
3.45292165529689	0.000883516926322344\\
3.47807955442474	0.000851416917682245\\
3.50323745355258	0.000822793458107796\\
3.53468482746239	0.000791026432371572\\
3.5661322013722	0.000762937075080802\\
3.59757957528201	0.000737896847193354\\
3.63531642397377	0.000711208684109673\\
3.67305327266554	0.000687597941717853\\
3.71079012135731	0.000666582255038549\\
3.75481644483104	0.000644843379388941\\
3.79884276830477	0.000625650467128816\\
3.84286909177851	0.000608629318184776\\
3.8931848900342	0.000591453091531196\\
3.94350068828989	0.00057636227361914\\
3.99381648654558	0.000563068155313394\\
4.05042175958324	0.000549970856654651\\
4.10702703262089	0.000538585085658983\\
4.16992178044051	0.000527682232182919\\
4.23281652826012	0.000518394717646721\\
4.29571127607974	0.000510538579622474\\
4.36489549868131	0.000503378529607062\\
4.43407972128289	0.000497627537678689\\
4.50326394388447	0.000493173932273565\\
4.57873764126801	0.000489700067952012\\
4.65421133865155	0.000487601975433486\\
4.72968503603508	0.000486838382582874\\
4.80515873341862	0.000487394934802753\\
4.88063243080216	0.000489282108793567\\
4.9561061281857	0.000492534351331454\\
5.03157982556924	0.000497210320941632\\
5.10076404817082	0.000502818744940555\\
5.16994827077239	0.000509779295870397\\
5.23913249337397	0.000518199387425746\\
5.30831671597555	0.000528211580887196\\
5.37121146379516	0.000538830516846672\\
5.43410621161478	0.000551043711295444\\
5.49700095943439	0.000565022091083254\\
5.55989570725401	0.000580963706520762\\
5.62279045507362	0.000599091327384098\\
5.67939572811128	0.000617469527136665\\
5.74229047593089	0.00064037791930343\\
5.80518522375051	0.000666005715225384\\
5.89323787069797	0.00070591099824999\\
5.96242209329955	0.000738290716733895\\
5.99386946720936	0.000751804107355189\\
6.0190273663372	0.000761140013686445\\
6.04418526546505	0.000768274269012833\\
6.06305368981093	0.00077148987506053\\
6.08192211415682	0.000772161728989422\\
6.09450106372074	0.000770796400098993\\
6.10708001328466	0.000767643626607293\\
6.11965896284859	0.000762365717079585\\
6.13223791241251	0.000754578666837666\\
6.14481686197643	0.000743852194374284\\
6.15739581154036	0.000729714442414746\\
6.16997476110428	0.000711665422001048\\
6.17626423588624	0.000701017156795248\\
6.1825537106682	0.00068920749696701\\
6.18884318545016	0.000676185296276543\\
6.19513266023212	0.000661911559276139\\
6.20142213501409	0.000646367162823111\\
6.21400108457801	0.000611568311908615\\
6.22658003414193	0.000572707800602471\\
6.24544845848782	0.000516001545057467\\
6.25173793326978	0.000504446353207446\\
6.25802740805174	0.000506426786009889\\
6.2643168828337	0.000545404807832735\\
6.27060635761566	0.000702138125818613\\
6.27689583239763	0.00108732675635415\\
6.28318530717959	0.00079311708260599\\
};
\addlegendentry{TF-IRKA}

\end{axis}
\end{tikzpicture}%

%% file: figures/E2Cplus.tex
%
%
\begin{tikzpicture}

\begin{axis}[%
  width=5cm,
  height=5cm,
  at={(0,0)},
scale only axis,
xmin=4,
xmax=18,
xtick={4, 6, 8, 10, 12, 14, 16, 18},
ymode=log,
ymin=1e-05,
ymax=0.1,
yminorticks=true,
axis background/.style={fill=white},
ylabel={$\etwod$ relative error},
legend style={at={(1,1)}, legend cell align=left, align=left, draw=white!15!black}
]
\addplot [color=blue, mark=*, mark options={solid, blue}, line width = 1.5pt]
  table[row sep=crcr]{%
4	0.02843034808139\\
6	0.00913032108222961\\
8	0.00258748018618611\\
10	0.000648000460536056\\
12	0.000240800413644906\\
14	0.000122141962027006\\
16	5.54071763337514e-05\\
18	1.67025511088851e-05\\
};
\addlegendentry{TF-IRKA}

\addplot [color=red, mark=square*, mark options={solid, red}, line width = 1.5pt]
  table[row sep=crcr]{%
4	0.0168502327508943\\
6	0.0128602867131481\\
8	0.00258745863643198\\
10	0.000647998614438113\\
12	0.000240800285542113\\
14	8.8413519643049e-05\\
16	6.03881239890363e-05\\
18	1.67025509973624e-05\\
};
\addlegendentry{Algorithm 1}

\end{axis}
\end{tikzpicture}%

%% file: figures/H2AerrorBDF2.tex
%
%
\definecolor{mycolor1}{rgb}{0.4940 0.1840 0.5560}%

\begin{tikzpicture}

\begin{axis}[%
width=5.5cm,
height=4cm,
at={(0,0)},
scale only axis,
xmin=4,
xmax=40,
xtick={4, 8, 12, 16, 20, 24, 28, 32, 36, 40},
ymode=log,
ymin = 1e-4,
ymax=1e-1,
yminorticks=true,
xlabel = {$r$},
ylabel = {\small $\durenhardy$ relative error},
axis background/.style={fill=white},
legend columns = 2,
legend style={at={(1,1.35)}, legend cell align=left, align=left, draw=white!15!black}
]
\addplot [color=blue, mark=*, mark options={solid, blue}, line width = 1.5pt]
  table[row sep=crcr]{%
4	0.0435393164548685\\
8	0.0292079770950595\\
12	0.0123649290178558\\
16	0.00673109368877329\\
20	0.00544937901251946\\
24	0.00330748769688569\\
28	0.00312674515910364\\
32	0.00150708977694266\\
36	0.00121822887526225\\
40	0.000324511199229888\\
};
\addlegendentry{\small AAA on $\unitdisk$}

\addplot [color=mycolor1, dashed, mark=o, mark size = 4pt, mark options={solid, mycolor1}, line width = 1.5pt]
  table[row sep=crcr]{%
4	0.0437971274316339\\
8	0.0299006244625715\\
12	0.0150299435670798\\
16	0.0107906561031024\\
20	0.00776838919783458\\
24	0.00617774091943019\\
28	0.00605055803743892\\
32	0.00371909561924997\\
36	0.00321356638925865\\
40	0.000463776107591807\\
};
\addlegendentry{\small AAA on $\A$}


\addplot [color=red, mark=square*, mark options={solid, red}, line width = 1.5pt]
  table[row sep=crcr]{%
4	0.0425427302065721\\
8	0.0339979322925562\\
12	0.0122266207420413\\
16	0.00233790968122407\\
20	0.00210437704365361\\
24	0.00141415596257468\\
28	0.000735672479244166\\
32	0.000464981688940679\\
36	0.000438932523808747\\
40	0.000196238788536356\\
};
\addlegendentry{\small Algorithm 1}

\end{axis}
\end{tikzpicture}%

%% file: figures/errorBDF.tex
%
%
\begin{tikzpicture}

\begin{axis}[%
width=6cm,
height=2cm,
at={(0,3.5cm)},
scale only axis,
xmin=0,
xmax=100,
ymin=0,
ymax=0.25,
xticklabel=\empty,
axis background/.style={fill=white},
legend style={at={(1,1)}, legend cell align=left, align=left, draw=white!15!black}
]
\addplot [color=black, line width = 1.5pt]
  table[row sep=crcr]{%
1	0\\
9	0\\
10	0.106643326466155\\
11	0.190239420372095\\
12	0.205684124072121\\
13	0.190055893413557\\
14	0.180856246957717\\
15	0.185931521690108\\
16	0.197469127503751\\
17	0.208923180709931\\
18	0.218939722712364\\
19	0.227948666836156\\
20	0.23537172864981\\
21	0.239565228628507\\
22	0.238970549841042\\
23	0.232964522295759\\
24	0.222015047859699\\
25	0.207331921629176\\
26	0.190316468645506\\
27	0.172122064902723\\
28	0.153536913571742\\
29	0.135143933596481\\
30	0.117499739479484\\
31	0.101127979152878\\
32	0.0863728894892688\\
33	0.0733154091014256\\
34	0.0618569808230234\\
35	0.0518738605495912\\
36	0.0432855792973044\\
37	0.0360085877697571\\
38	0.0298985913224357\\
39	0.0247684755000819\\
40	0.0204544864383394\\
41	0.0168480480772359\\
42	0.0138673520735892\\
43	0.0114187570006266\\
44	0.00939733797308406\\
45	0.0077164958051128\\
46	0.0063232082146385\\
47	0.00518242531184399\\
48	0.00425489121487033\\
49	0.00349492220657055\\
50	0.00286463895790234\\
51	0.00234248425084616\\
52	0.00191620739161635\\
54	0.00129039402142439\\
56	0.000863764199635852\\
58	0.000578737983730093\\
61	0.000317929439944464\\
65	0.000143193029316535\\
71	4.29907447596634e-05\\
85	2.6131706505339e-06\\
100	1.29896761791315e-07\\
};
\addlegendentry{\small $y_k$}

\addplot [color=red, dashed, line width = 1.5pt]
  table[row sep=crcr]{%
1	0\\
9	0\\
10	0.106643316717026\\
11	0.190239506915518\\
12	0.205683850247979\\
13	0.190056224413979\\
14	0.18085654110611\\
15	0.185929814170152\\
16	0.197471571552626\\
17	0.208922501492054\\
18	0.218937719415734\\
19	0.227950515885823\\
20	0.235373007433751\\
21	0.239562759529278\\
22	0.238969258443987\\
23	0.232966934774325\\
24	0.222016671282972\\
25	0.207330499491377\\
26	0.190315917844558\\
27	0.172124118483495\\
28	0.153536997218552\\
29	0.135139250851225\\
30	0.117495203726477\\
31	0.101129537740277\\
32	0.0863792174747857\\
33	0.0733207270750142\\
34	0.0618589796725928\\
35	0.0518739794040215\\
36	0.0432834221842739\\
37	0.0360018469243357\\
38	0.0298891207482654\\
39	0.024763415917576\\
40	0.0204584780335608\\
41	0.0168570448724523\\
42	0.0138732877518777\\
43	0.0114190485058288\\
44	0.00939664443566812\\
45	0.00771925119650518\\
46	0.00632764812750963\\
47	0.00518356017305166\\
48	0.00425128972196376\\
49	0.00349003358080324\\
50	0.00286185336248934\\
51	0.00234126270537161\\
52	0.00191395657915905\\
54	0.00128710618260186\\
56	0.000864303517246867\\
58	0.000578606552807059\\
61	0.000320674925447406\\
65	0.000145948250320771\\
71	4.50037495625111e-05\\
84	3.47877954709475e-06\\
100	1.53515358647383e-07\\
};
\addlegendentry{\small $\widehat{y}_k$}
\end{axis}

\begin{axis}[%
width=6cm,
height=3cm,
at={(0,0)},
scale only axis,
xmin=0,
xmax=100,
ymode=log,
ymin=1e-9,
ymax=0.001,
yminorticks=true,
axis background/.style={fill=white},
xlabel = {$k$},
legend columns = 2,
legend style={at={(0.98,0.25)}, legend cell align=left, align=left, draw=white!15!black}
]
\addplot [color=red, dashed, line width = 1.5pt]
  table[row sep=crcr]{%
10	7.74241661775712e-08\\
11	5.45878436349414e-07\\
12	1.88499467166321e-06\\
13	4.41495542592152e-06\\
14	8.10017549118487e-06\\
15	1.10502400788605e-05\\
16	9.07008351995695e-06\\
17	7.07932431182775e-06\\
18	1.0801511375141e-05\\
19	6.67774825172859e-06\\
20	1.09280187941202e-05\\
21	9.21100526872087e-06\\
22	1.17880668067941e-05\\
23	9.83422340779182e-06\\
24	1.06507672777236e-05\\
25	5.52509948291219e-06\\
26	5.70428615899496e-06\\
27	9.11950574871995e-06\\
28	1.35098200938377e-05\\
29	2.03539836434848e-05\\
30	2.54358610882788e-05\\
31	2.31181726625946e-05\\
32	3.01051632846415e-05\\
33	2.41423190300293e-05\\
34	1.60944592206123e-05\\
35	1.71264869632892e-05\\
36	2.05841413831286e-05\\
37	3.07130894949121e-05\\
38	4.05800242391201e-05\\
39	3.21182273041244e-05\\
40	2.94339821102342e-05\\
41	3.80682473560952e-05\\
42	2.55972277323352e-05\\
43	4.90120205851929e-06\\
44	1.41122626943096e-05\\
45	2.09525448639984e-05\\
46	1.87044415628863e-05\\
47	2.10113965856973e-05\\
48	2.54380031621614e-05\\
49	1.90819320012069e-05\\
50	1.15145797236815e-05\\
51	4.96083256967405e-06\\
52	1.56380092692606e-05\\
53	2.03698907245235e-05\\
54	1.17077963937472e-05\\
55	8.5405206892066e-06\\
56	7.94181567811304e-06\\
57	6.12032591144337e-07\\
58	2.52230622806916e-06\\
59	4.52201915401775e-06\\
60	1.31825924494758e-05\\
61	1.4490805294999e-05\\
62	9.60631454576609e-06\\
63	7.1716408924713e-06\\
64	9.97051341431236e-06\\
65	1.38011646539041e-05\\
66	1.35492382133166e-05\\
67	9.69213830519268e-06\\
68	7.33306850899053e-06\\
69	8.3250848488846e-06\\
70	9.83773544722713e-06\\
71	9.14483908982564e-06\\
72	6.52115814338654e-06\\
73	4.79871133751622e-06\\
74	5.1248192315693e-06\\
75	5.73472004969452e-06\\
76	5.20019276531778e-06\\
78	2.57153592936402e-06\\
79	2.69408875826352e-06\\
80	2.97827651386468e-06\\
81	2.68406363548551e-06\\
82	1.85573729108882e-06\\
83	1.22122281869562e-06\\
84	1.28235343827545e-06\\
85	1.4426880893371e-06\\
86	1.31118685473747e-06\\
87	8.9501882715988e-07\\
88	5.30889828223636e-07\\
89	5.67954021055676e-07\\
90	6.70802219742772e-07\\
91	6.22569670353056e-07\\
92	4.22680148041125e-07\\
93	2.14694737414809e-07\\
94	2.36909645832812e-07\\
95	3.05539756430441e-07\\
96	2.92574772759138e-07\\
97	2.00112471709279e-07\\
98	8.27905238384242e-08\\
99	9.33444672593281e-08\\
100	1.38506120792309e-07\\
};
\addlegendentry{\small $|y_k-\widehat{y}_k|$}
\addplot [color=black, line width = 1.5pt]
  table[row sep=crcr]{%
1	0.000131257205685338\\
100	0.000131257205685338\\
};
\addlegendentry{\small $\|G-G_r\|_{\durenhardy}$}
\end{axis}
\end{tikzpicture}%

%% file: figures/H2errorBDF.tex
%
%
\begin{tikzpicture}

\begin{axis}[%
width=4cm,
height=5.5cm,
at={(1.011in,0.642in)},
scale only axis,
xmin=4,
xmax=14,
ymode=log,
ymin=1e-07,
ymax=0.02,
yminorticks=true,
axis background/.style={fill=white},
xlabel = {$r$},
ylabel = {$\durenhardy$ relative error},
]
\addplot [color=red, mark=square*, line width = 1.5pt, mark options={solid, red}]
  table[row sep=crcr]{%
4	0.0108446616494035\\
6	0.0021019330977888\\
8	0.000251359290588723\\
10	3.38690806299531e-05\\
12	4.05595332015679e-06\\
14	5.23591491065137e-07\\
};

\end{axis}
\end{tikzpicture}%